\pdfoutput=1 

\newif\REVIEW

\ifdefined\REVIEW
    \documentclass[
        a4paper,
        onecolumn
    ]{article}
    \usepackage{lineno, setspace}
\else
    \documentclass[
        a4paper,
        twocolumn
    ]{article}
\fi

\usepackage{IEK10} 
\usepackage{natbib}

\usepackage{upgreek}
\usepackage{bbm}

\DeclareSIUnit\year{yr}
\DeclareSIUnit\annum{a}
\DeclareSIUnit\ton{t}
\DeclareSIUnit{\million}{\text{M}}
\newcommand{\Edot}{\dot{E}}

\def\IEK10{
  Institute of Climate and Energy Systems,
  Energy Systems Engineering (ICE-1),
  Forschungszentrum J\"ulich GmbH,
  J\"ulich 52425,
  Germany
}
\def\RWTH{
  RWTH Aachen University,
  Aachen 52062,
  Germany
}
\def\JARA{
  JARA-ENERGY,
  J{\"u}lich 52425,
  Germany
}
\def\SVT{
  RWTH Aachen University,
  Process Systems Engineering (AVT.SVT),
  Aachen 52074,
  Germany
}

\newcommand{\mytitle}{Robust Design of Multi-Energy Systems Accounting for Mixed-Integer Operational Problems}

\newcommand{\affil}{
  \begin{itemize}[leftmargin=3mm, itemsep=0mm]
    \item[$^a$]\IEK10
    \item[$^b$]\RWTH
    \item[$^c$]\SVT
    \item[$^d$]\JARA
  \end{itemize}
}

\def\firstAuthor{Moritz Wedemeyer}
\newcommand{\myauthor}{\firstAuthor$^{a,b}$,
Alexander Mitsos$^{d,a,c}$,
Manuel Dahmen$^{a,*}$
}
\newcommand{\correspondingauthor}{M. Dahmen, \IEK10 \newline E-mail: m.dahmen@fz-juelich.de}
\author{\myauthor}

\usepackage[
  colorlinks,
  linkcolor=blue,
  citecolor=blue,
  urlcolor=blue,
  pdftitle={\mytitle},
  pdfauthor={\firstAuthor}
]{hyperref}
\usepackage[capitalise, nameinlink]{cleveref}
\crefname{table}{Tab.}{Tab.}

\usepackage[colorinlistoftodos]{todonotes}

\begin{document}
\twocolumn[
\begin{@twocolumnfalse}

  \thispagestyle{firststyle}

  \begin{center}
    \begin{large}
      \textbf{\mytitle}
    \end{large} \\
    \myauthor
  \end{center}

  \vspace{0.5cm}

  \begin{footnotesize}
    \affil
  \end{footnotesize}

  \vspace{0.5cm}
Identifying robust designs for multi-energy systems is computationally challenging.
As rigorous approaches are often computationally intractable, heuristics are employed to generate candidate designs.
Specifically, we consider a heuristic that iteratively identifies and adds extreme scenarios to the design problem.
We theoretically investigate how three common nonconvexities, i.e., piecewise-linear energy inflow-outflow relationships, minimum part-loads, and storage complementarity, affect the robustness of designs identified by this heuristic.
We find that, if surplus energy cannot be curtailed, any of these nonconvexities may cause the heuristic to fail.
If curtailment is allowed, storage complementarity does not compromise robustness, and convex piecewise-linear inflow-outflow relationships can be reformulated linearly.
However, minimum part-loads may lead to failure of the heuristic.
Furthermore, if the optimal value function of the operational problem is nonconvex in the uncertain variables, the heuristic may fail.
We demonstrate these findings using an illustrative multi-energy system case study, in which minimum part-loads and nonconvex dependence of the objective function on heat-pump efficiency are identified as possible failure modes.
We rigorously verify the robustness of a design identified using the heuristic and find a scenario where, in one time step, \SI{1.1}{\kilo\watt} of the \SI{313.6}{\kilo\watt} electricity demand could not be satisfied.
However, when the number of representative scenarios is increased from $4$ to $6$ or $8$, the resulting designs are robust.
This demonstrates that the heuristic can serve as an effective first step in identifying robust designs.
However, when robustness guarantees are required, a rigorous solution method must be employed.

\vspace*{5mm}


\end{@twocolumnfalse}
]
\ifdefined\REVIEW
  \onecolumn
\else
  \newpage
\fi

\section{Introduction}
When designing multi-energy systems, it is paramount that these are robust, i.e., they can meet any energy demand that the system is expected to encounter during operation.
Robustness is an even stronger concern in isolated energy systems, which are not connected to a superordinate grid.
The limited number of users and small geographic extent of isolated systems lead to large relative fluctuations in load and renewable generation, necessitating robustness against uncertainties in both demand and renewable supply.
To hedge against uncertainty, multi-energy systems are often designed using historical realizations of the uncertainties.
As incorporating the full historical dataset is often computationally intractable, time-series aggregation methods are frequently applied \citep{kotzurimpactdifferenttime2018, teichgraeberclusteringmethodsfind2019, fleschutzHiddenCostUsing2025}.

A major limitation of these aggregation approaches is that they often neglect extreme events, even though such events can strongly influence system costs \citep{bahlTimeseriesAggregationSynthesis2017}.
Consequently, several approaches have been developed to better account for extreme scenarios.
One option is to manually include extreme periods in the representative scenarios used for system design \citep{dominguez-munozSelectionTypicalDemand2011, vollAutomatedSuperstructurebasedSynthesis2013}.
In a similar vein, \cite{gabrielliRobustOptimalDesign2019} created synthetic worst-case scenarios with artificially increased demands and analyzed the robustness of the candidate designs.
However, the interactions between different uncertainties are difficult to assess a priori.

In a series of papers, Bardow and co-workers \citep{bahlTimeseriesAggregationSynthesis2017, bahlRigorousSynthesisEnergy2018, baumgartnerRiSES3RigorousSynthesis2019, baumgartnerRiSES4RigorousSynthesis2019} developed methods that iteratively refine the temporal resolution of aggregated time series to accurately approximate the optimal objective value of the full dataset for mixed-integer linear programming (MILP) energy system design problems.
\cite{wangRiNSES4RigorousNonlinear2024} extended these methods to nonlinear problems.
In addition to improving the approximation of system costs, these approaches also aim to ensure the feasibility of the candidate design for the complete historical dataset.
In a similar vein, \cite{teichgraeberExtremeEventsTime2020} introduced the feasibility time-step heuristic, based on work by \cite{bahlRigorousSynthesisEnergy2017}.
Instead of increasing the resolution of the entire time series, as done in \cite{bahlTimeseriesAggregationSynthesis2017, bahlRigorousSynthesisEnergy2018} and \cite{baumgartnerRiSES3RigorousSynthesis2019, baumgartnerRiSES4RigorousSynthesis2019}, \cite{teichgraeberExtremeEventsTime2020} and \cite{bahlRigorousSynthesisEnergy2017} iteratively check the operational problem for feasibility and add violating scenarios to the design problem.

However, feasibility with respect to the historical data alone may not be sufficient to guarantee feasibility for future uncertainty realizations \citep{wedemeyerRobustEnergySystem2025}, since the available historical data represent only a finite subset of all possible operating conditions.
To obtain robustness guarantees, robust optimization \citep{ben-talRobustOptimization2009}, originally introduced by \cite{soysterTechnicalNoteConvex1973}, can be employed.
Robust optimization addresses optimization problems with uncertain parameters, e.g., future energy demands, renewable generation, and energy prices.
The objective is to identify solutions that remain feasible for all realizations of the uncertain parameters within a predefined uncertainty set.
Enforcing feasibility for all possible uncertainty realizations naturally leads to semi-infinite programs (SIPs) \citep{steinBilevelStrategiesSemiinfinite2003, charnesDualityHaarPrograms1962}, which involve finitely many decision variables but infinitely many constraints.

In many energy system applications, some decisions must be taken before the uncertainty is realized, while others can be adjusted afterwards.
This setting is commonly addressed using adjustable robust optimization \citep{ben-talAdjustableRobustSolutions2004}, in which here-and-now decisions are complemented by wait-and-see variables acting as recourse decisions.
Including recourse decisions leads to multilevel optimization problems, making adjustable robust optimization problems even harder to solve than regular robust optimization problems.

Existing approaches, therefore, typically balance the tractability of the robust optimization problem against the generality of the uncertainty set and the flexibility of the recourse decisions.
For example, \cite{bertsimasFiniteAdaptabilityMultistage2010} introduced the field of finite adaptability, where instead of allowing arbitrary recourse decisions, a predefined number of contingency plans, i.e., fixed recourse decisions, can be made ahead of time.
After the uncertainty has been realized, the best recourse decisions out of this set can then be chosen.
Their framework provides a tradeoff between the fixed recourse of regular robust optimization and the computationally challenging arbitrary recourse in adjustable robust optimization.
\cite{majewskiTRusTTwostageRobustness2017} proposed a robust energy system design approach, but restricted the uncertainty set to intervals around nominal parameter values.
They reformulated the resulting SIP into a tractable single-level problem by considering only the extreme values of these intervals.
However, \cite{grossmannoptimumdesignchemical1978} have shown that this reformulation is not valid for general polyhedral uncertainty sets.
Similarly, \cite{abdinOptimizingGenerationExpansion2022} propose an adjustable robust optimization approach based on the same interval-type uncertainty set, in which wait-and-see variables are approximated using affine decision rules parametrized by here-and-now variables.
In contrast, \cite{shenLargescaleIndustrialEnergy2020} introduce a data-driven uncertainty set to better capture the correlation between uncertain parameters, but do not explicitly model recourse decisions.

Recently, we transferred a robust design approach from process systems engineering \citep{grossmannoptimumdesignchemical1978} to energy system design \citep{wedemeyerRobustEnergySystem2025}.
While this approach enables the exact solution of robust energy system design problems without restrictions on the structure of the uncertainty set or the recourse decisions, it is computationally demanding.
Consequently, in practical settings, heuristics are needed, such as the feasibility time-step heuristic introduced by \cite{teichgraeberExtremeEventsTime2020}.

We previously proved that for convex operational problems, the feasibility time-step heuristic indeed guarantees feasibility for all scenarios within the convex hull of the historical data \citep{wedemeyerRobustEnergySystem2025}. 
However, multi-energy system design problems are often formulated as MILP problems, which are inherently nonconvex \citep{ommenComparisonLinearMixed2014,kotzurModelersGuideHandle,wirtzDesignOptimizationMultienergy2021,hoffmannReviewMixedintegerLinear2024, mancoReviewMultiEnergy2024}.
As a result, the robustness guarantees established for convex operational problems no longer apply.
Nevertheless, the original work by \cite{bahlRigorousSynthesisEnergy2017} already applied the heuristic to an MILP energy system synthesis problem.

This raises the question of whether the feasibility time-step heuristic can provide robustness guarantees in the presence of nonconvex operational constraints.
In our previous work \citep{wedemeyerRobustEnergySystem2025}, we showed that minimum part-load operation, if combined with a no-curtailment assumption, can cause the feasibility time-step heuristic \citep{teichgraeberExtremeEventsTime2020} to identify designs that are not robust, i.e., designs for which not all realizations in the uncertainty set are feasible.
This minimal example demonstrated one possible failure mode of the heuristic.
However, it remained unclear whether other practically relevant nonconvexities can cause similar failures, whether some nonconvexities are harmless, and whether these failures persist when curtailment is allowed.

The present manuscript addresses this gap by systematically analyzing the conditions under which the feasibility time-step heuristic guarantees robustness despite the presence of three common nonconvexities in multi-energy system models, i.e., nonconvex minimum part-load constraints, piecewise-linear performance curves, and storage complementarity constraints.
Specifically, we analyze the effects of these constraints on the robustness guarantees of the heuristic under both curtailment and no-curtailment assumptions.
We extend our previous work \citep{wedemeyerRobustEnergySystem2025} by providing additional minimal examples showing that failures of the feasibility time-step heuristic may also occur when curtailment is allowed.
We then go beyond counterexamples by theoretically analyzing which of the investigated nonconvexities invalidate the robustness guarantee and which preserve it under specific conditions.
In addition, we propose a rigorous verification method to test whether candidate designs obtained with the heuristic are robust.

Finally, an illustrative case study demonstrates the practical relevance of the theoretical analysis.
It shows how one of the theoretically identified failure modes can arise in a multi-energy system design problem and cause the feasibility time-step heuristic to identify a non-robust design.

The remainder of this work is structured as follows.
Section \ref{sec:PS} presents the modeling framework for multi-energy system design and discusses three key nonconvexities arising in the operational problem.
Their impact on robustness is then analyzed.
Section \ref{sec:case_study} introduces an illustrative case study that demonstrates the impact of nonconvexities.
Additionally, a method is proposed to verify whether a design candidate is robust for all scenarios within the convex hull of the historical data, and its application is demonstrated using the case study.
Finally, Section \ref{sec:conclusion} discusses the results and summarizes the main findings.
\section{Problem Structure}\label{sec:PS}
Throughout the manuscript, scalar-valued quantities are denoted in regular font, e.g., $x$, vector-valued quantities are denoted in bold font, e.g., $\mathbf{x}$, and set-valued quantities are denoted in calligraphic font, e.g., $\mathcal{X}$.
We consider energy system design problems of the form:
\begin{align*}\tag{PS}\label{problem_structure}
    \underset{\mathbf{x} \in \mathcal{X}, \mathbf{z}_{s} \in \mathcal{Z}_{s}(\mathbf{x}) \forall s \in \mathcal{S}}{\min} & \quad c_{inv}(\mathbf{x}) + \sum_{s \in \mathcal{S}} c_{op}(\mathbf{z}_{s}) &\\
    \text{s.t.} & \quad \underset{\mathbf{y} \in \mathcal{Y}}{\max}\underset{\mathbf{z} \in \mathcal{Z}(\mathbf{x}, \mathbf{y})}{\min} \quad E_{gap}(\mathbf{x}, \mathbf{y}, \mathbf{z}) \le 0&\\
\end{align*}
with feasible sets
\begin{align*}
    \mathcal{X} = \{\mathbf{x} \in \mathbb{R}^{n_x} | \ \mathbf{g}_{x}(\mathbf{x}) \le \mathbf{0} \ (Design \ constraints)\}\\
    \mathcal{Z}_{s}(\mathbf{x}) = \{\mathbf{z}_{s} \in \mathbb{R}^{n_z} | \ \mathbf{g}_{en}(\mathbf{x}, \mathbf{y}_{s}, \mathbf{z}_{s}) \le \mathbf{0} \ (Energy \ system \ model)\}\\
    \mathcal{Y} = \{\mathbf{y} \in \mathbb{R}^{n_y} | \ \mathbf{g}_{y}(\mathbf{y}) \le \mathbf{0} \ (Uncertainty \ bounds)\}\\
    \mathcal{Z}(\mathbf{x}, \mathbf{y}) = \{\mathbf{z} \in \mathbb{R}^{n_z} | \ \mathbf{g}_{en}(\mathbf{x}, \mathbf{y}, \mathbf{z}) \le \mathbf{0} \ (Energy \ system \ model)\},
\end{align*}
where $\mathbf{x}$ represents design decisions such as component capacities, $\mathbf{y}_{s}$ are the fixed uncertainty realizations of the representative scenarios $s$, and $\mathbf{z}_{s}$ represents the corresponding upper-level operational decisions such as component energy outflows.
$c_{inv}$ and $c_{op}$ are the investment and operational costs, respectively.
$\mathbf{g}_{x}(\mathbf{x})$ are the design constraints, e.g., piecewise-linear investment curves, and $\mathbf{g}_{y}(\mathbf{y})$ are the uncertainty bounds, i.e., constraints describing the convex hull of the historical data.
$\mathbf{g}_{en}(\mathbf{x}, \cdot, \cdot) \le \mathbf{0}$ are the model equations that describe the behavior of the energy system.
\begin{align*}\tag{MLP}\label{medial_level}
     \underset{\mathbf{y} \in \mathcal{Y}}{\max}\underset{\mathbf{z} \in \mathcal{Z}(\mathbf{x}, \mathbf{y})}{\min} \quad E_{gap} (\mathbf{x}, \mathbf{y}, \mathbf{z})
\end{align*}
The embedded optimization problem, referred to as the medial-level problem \eqref{medial_level}, verifies whether the candidate design remains feasible for all uncertainty realizations $\mathbf{y}$, such as realizations of wind power, by allowing operational recourse decisions $\mathbf{z}$ to adapt system operation to each realization.
We refer to the objective function of problem \eqref{medial_level} as the energy gap, which is defined as:
\begin{equation}\label{eq:mlp_obj}
    \begin{aligned}
    &E_{gap} (\mathbf{x}, \mathbf{y}, \mathbf{z}) = \\
    &\underset{ e \in \mathcal{E}, t \in \mathcal{T}}{\max} \left\{ \Edot_{e, dem, t} - \sum_{c \in \mathcal{C}_{e, out}} \Edot_{e, c, t} + \sum_{c \in \mathcal{C}_{e, in}} \Edot_{e, c, t} + \sum_{c \in \mathcal{C}_{e, sto}} \left(\Edot_{e, c, t, in} - \Edot_{e, c, t, out} \right)\right\}
    \end{aligned}
\end{equation}
Here, the energy gap represents the maximum energy balance violation across all energy forms $e \in \mathcal{E}$ and time steps $t \in \mathcal{T}$.
$\Edot_{e, dem, t}$ is the energy demand for energy form $e$ at time $t$, and $\Edot_{e, c, t}$ represent energy flows associated with component 
$c$, where $\mathcal{C}_{e, in}$ and $\mathcal{C}_{e, out}$ denote the sets of components consuming and supplying energy form $e$, respectively.
For storage components $c \in \mathcal{C}_{e, sto}$, $\Edot_{e, c, t, in}$ and $\Edot_{e, c, t, out}$ denote the energy inflows and outflows.

\eqref{problem_structure} is an adjustable robust optimization problem \citep{ben-talAdjustableRobustSolutions2004} which can naturally be formulated as a so-called existence-constrained semi-infinite program (ESIP) \citep{djelassiGlobalSolutionSemiinfinite2021}.
Importantly, if the feasible set of the lower-level problem $\mathcal{Z}(\mathbf{x}, \mathbf{y})$ depends on the uncertain variables $\mathbf{y}$, \eqref{medial_level} is a generalized semi-infinite program (GSIP).
To facilitate the solution of \eqref{medial_level}, the coupling constraints $g_{j}(\mathbf{y}, \mathbf{z})$ in $\mathcal{Z}(\mathbf{x}, \mathbf{y})$, which introduce the dependence on the medial-level variables $\mathbf{y}$, are moved into the objective function:
\begin{align}\tag{MLP-ref}\label{eq:mlp_ref}
     \underset{\mathbf{y} \in \mathcal{Y}}{\max}\underset{\mathbf{z} \in \mathcal{Z}(\mathbf{x})}{\min} \quad \max \left\{E_{gap} (\mathbf{x}, \mathbf{y}, \mathbf{z}), \underset{j \in \mathcal{J}}{\max}\left\{g_{j}(\mathbf{y}, \mathbf{z})\right\}\right\}
\end{align}
Here $j \in \mathcal{J}$ are the indices of the coupling constraints.
If \eqref{medial_level} is replaced  by \eqref{eq:mlp_ref} in \eqref{problem_structure}, an equivalent ESIP is obtained (cf. Section $5.2.4$ in \cite{djelassidiskretisierungsbasiertealgorithmenfur2020}).

If $\mathcal{Y}$ is a polyhedral set, as we will assume here, a sufficient condition for the worst-case scenario to lie on a vertex of the uncertainty set is convexity of the optimal value function
\begin{equation*}
    \nu^{*}(\mathbf{y}) = \underset{\mathbf{z} \in \mathcal{Z}(\mathbf{x})}{\min} \quad \max \left\{E_{gap} (\mathbf{x}, \mathbf{y}, \mathbf{z}), \underset{j \in \mathcal{J}}{\max}\left\{g_{j}(\mathbf{y}, \mathbf{z})\right\}\right\}
\end{equation*}
with respect to $\mathbf{y}$.
Sufficient conditions for convexity of $\nu^{*}(\mathbf{y})$ are joint convexity of the objective function $\max \left\{E_{gap} (\mathbf{x}, \mathbf{y}, \mathbf{z}), \underset{j \in \mathcal{J}}{\max}\left\{g_{j}(\mathbf{y}, \mathbf{z})\right\}\right\}$ in $\mathbf{y}$ and $\mathbf{z}$, and convexity of $\mathcal{Z}(\mathbf{x})$ (cf. Section $3.1.5$ in \cite{boydConvexOptimization2004}).

Consequently, if $\nu^{*}(\mathbf{y})$ is nonconvex, the worst-case scenario may not lie on a vertex of the uncertainty set, which may lead to failure of the feasibility time-step heuristic and thus a non-robust design \citep{wedemeyerRobustEnergySystem2025}.

\eqref{problem_structure} can be solved using adaptive discretization approaches \citep{djelassidiskretisierungsbasiertealgorithmenfur2020} based on the algorithm by \cite{blankenshipInfinitelyConstrainedOptimization1976}, which is similar to the more recently proposed column-and-constraint generation approach introduced by \cite{zengSolvingTwostageRobust2013} and has also been adapted for the solution of mixed-integer nonlinear bilevel programs \citep{mitsosGlobalSolutionNonlinear2010}.
\begin{figure}
    \centering    \includegraphics[width=0.7\linewidth]{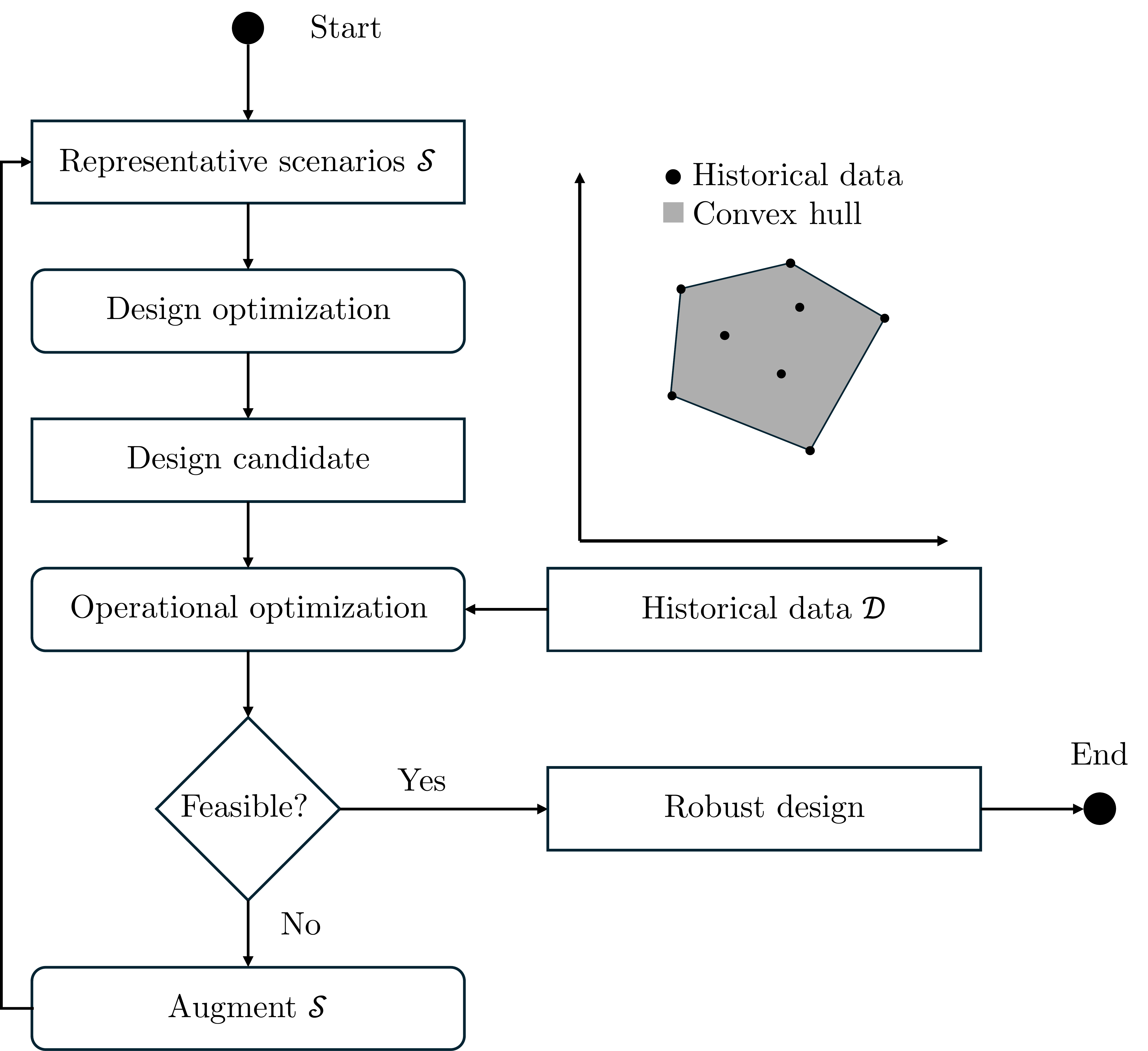}
    \caption{Illustration of the feasibility time-step heuristic \citep{teichgraeberExtremeEventsTime2020}: A candidate design is obtained using a set of representative scenarios $\mathcal{S}$. Its operational feasibility is evaluated for all historical scenarios (circles) in the historical data $\mathcal{D}$. If all scenarios are feasible, the design is deemed robust; otherwise, the set of representative scenarios is augmented with one or more infeasible scenarios from the historical data, and the procedure is repeated. For convex operational problems, the feasibility of all historical scenarios guarantees feasibility for the entire interior of their convex hull (shaded region) \citep{wedemeyerRobustEnergySystem2025}.}\label{fig:feasts}
\end{figure}
As shown in our previous work \citep{wedemeyerRobustEnergySystem2025}, these adaptive discretization approaches are based on a principle similar to the feasibility time-step heuristic, which is illustrated in Figure \ref{fig:feasts}.
In both methods, design optimization is performed using a set of representative scenarios $\mathcal{S}$ drawn from the historical dataset $\mathcal{D}$ to obtain candidate designs.
Subsequently, worst-case scenarios are iteratively identified and added to $\mathcal{S}$ to refine the design problem.

The two approaches differ in how the worst-case scenarios are determined.
In the feasibility time-step heuristic, operational feasibility is evaluated for all historical scenarios.
If one or more infeasible scenarios are identified, they are added to the representative scenario set, and a new candidate design is computed (cf. Figure \ref{fig:feasts}).
The procedure terminates once all historical scenarios are feasible, at which point the design is deemed robust with respect to the historical data.

In contrast, adaptive discretization algorithms \citep{djelassidiskretisierungsbasiertealgorithmenfur2020} identify worst-case scenarios by solving \eqref{eq:mlp_ref} directly.
Consequently, if the uncertainty set $\mathcal{Y}$ is defined as the convex hull of the historical data, the feasibility time-step heuristic solves a relaxation of \eqref{problem_structure} \citep{wedemeyerRobustEnergySystem2025}.
\subsection{Common nonconvexities}\label{sec:ncvx}
Based on recent publications and reviews of multi-energy system modeling \citep{ommenComparisonLinearMixed2014, kotzurModelersGuideHandle, wirtzDesignOptimizationMultienergy2021, hoffmannReviewMixedintegerLinear2024, mancoReviewMultiEnergy2024}, we have identified $5$ types to which nonconvexities occurring in energy system design problems can typically be attributed:
\begin{enumerate}
    \item Piecewise-linear investment curves
    \item Multiple units of the same component type
    \item Piecewise-linear performance curves
    \item Minimum part-loads
    \item Complementarity of storage input and output
\end{enumerate}

Note that the present work focuses on local multi-energy systems and on nonconvexities that arise within individual time steps.
Other sources of nonconvexity may occur in broader modeling settings.
For example, nonconvexities arise from AC power flow equations \citep{bienstockMathematicalProgrammingFormulations2020}.
However, the present study does not model energy transport and instead adopts a copper-plate assumption, under which energy can be exchanged freely within the system.
In larger-scale systems, intertemporal operational constraints that couple binary on/off variables across time steps can introduce additional nonconvexities \citep{wirtzDesignOptimizationMultienergy2021}.
Examples include limits on the number of start-ups and shut-downs, minimum up- and down-time constraints, and ramping constraints \citep{mancoReviewMultiEnergy2024}.
These additional sources of nonconvexity are beyond the scope of the present study.

Furthermore, we focus on the impact of nonconvexities on the robustness of energy system designs.
Piecewise-linear investment curves and integer variables used to model multiple units of the same component type appear only in the upper-level design problem \eqref{problem_structure} and not in the embedded problem \eqref{medial_level}, which ensures the robustness of the design.
Consequently, these nonconvexities do not affect the robustness of designs identified by the feasibility time-step heuristic.

This leaves $3$ relevant nonconvexities in the operational problem.
For these nonconvexities, their impact depends on two additional factors: (i) whether the input energy flow of the considered component is externally supplied and unconstrained, and (ii) whether curtailment of excess energy is permitted.

Multi-energy systems comprise components that convert between different energy carriers and must satisfy demand constraints for each carrier.
If the input energy form is externally supplied and unconstrained, and thus does not appear in any system balance equation, the nonconvexities affect only operational costs but not robustness, as the external supply is assumed to be able to meet any required demand.
If the externally supplied energy carrier is subject to any constraints, for example, limits on supply, it may affect robustness, and these constraints have to be included in the model.

Curtailment refers to the ability to dispose of surplus energy, for example, by releasing excess heat to the environment via a heat exchanger.
Curtailment has a substantial influence on the mathematical structure and difficulty of the problem.
If curtailment is allowed, \eqref{problem_structure} constitutes a regular ESIP \citep{djelassiGlobalSolutionSemiinfinite2021}, in which the objective function \eqref{eq:mlp_obj} of the medial-level problem \eqref{medial_level} must be less than or equal to zero.
If curtailment is prohibited, the energy balance must be satisfied exactly, resulting in a semi-infinite equality constraint.
However, the presence of a semi-infinite equality constraint violates the assumption that a Slater point exists \citep{djelassiGlobalSolutionSemiinfinite2021}.
The existence of Slater points, even near-optimal ones, is rather standard among deterministic global solution methods for semi-infinite problems, including for~\eqref{problem_structure}, as they rely on restriction to generate feasible bounds (upper bounding procedure).
In the absence of Slater points, the algorithms, including those proposed in \cite{djelassiGlobalSolutionSemiinfinite2021}, are not guaranteed to converge, and typically do not do so in practice.
In the presence of semi-infinite equality constraints, specialized solution strategies are required, e.g., \cite{djelassiDiscretizationbasedAlgorithmsGeneralized2019} or \cite{stuberSemiInfiniteOptimizationImplicit2015}.
In line with previous works on the robustness of energy systems \citep{majewskiTRusTTwostageRobustness2017}, we focus on the case with curtailment, as it yields a more tractable yet realistic optimization problem formulation.
In practice, excess energy can often be dissipated, e.g., via curtailing renewable energy sources or dissipating heat energy to the environment through the use of heat exchangers.

In addition to the previously mentioned nonconvexities in the operational problem, the way uncertainties enter the operational problem can also lead to the feasibility time-step heuristic being unable to identify a robust design.
Typical uncertainties in energy system optimization problems include energy demand and weather realizations.

While energy demand enters the operational problem linearly in the energy balance equations of the objective function \eqref{eq:mlp_obj}, weather variables may appear nonlinearly.
For example, the efficiency of a heat pump $\eta_{hp, nom, t}(T_{amb, t})$ may depend on ambient temperature $T_{amb, t}$ \citep{sassModelCompendiumData2020} and appears in the energy inflow-outflow relation as:
\begin{equation*}
    \Edot_{e_{in}, hp, t} - \frac{\Edot_{e_{out}, hp, t}}{\eta_{hp, nom, t}(T_{amb, t})} = 0
\end{equation*}
When $\Edot_{e_{in}, hp, t}$ is replaced in the energy balance by the expression $\frac{\Edot_{e_{out}, hp, t}}{\eta_{hp, nom, t}(T_{amb, t})}$, a nonconvex dependence of the objective function \eqref{eq:mlp_obj} on the uncertain variable $\eta_{hp, nom, t}(T_{amb, t})$ arises.
Specifically, the efficiency of the heat pump is defined as $\eta_{hp, nom, t}(T_{amb, t}) = \frac{0.36T_{hp}}{T_{hp} - T_{amb, t}}$, where $T_{hp}$ is the output temperature of the heat pump.
This leads to the nonconvex term $\frac{\Edot_{e_{out}, hp, t} \left( T_{hp} - T_{amb, t}\right)}{0.36 T_{hp}}$.
Since $T_{hp}$ is a parameter, this expression can be separated into a linear term $\frac{\Edot_{e_{out}, hp, t} T_{hp}}{0.36 T_{hp}}$ and the nonconvex bilinear term $ \frac{- \Edot_{e_{out}, hp, t} T_{amb, t} }{0.36 T_{hp}}$.
The latter term couples the operational variable $\Edot_{e_{out},hp,t}$ with the uncertain variable $T_{amb,t}$, and the Hessian with respect to these two variables is indefinite for all values of the variables.
Hence, the term is not jointly convex on any open domain in which both variables can vary.
Therefore, the energy gap \eqref{eq:mlp_obj} contains a term that is not jointly convex in $\mathbf{y}$ and $\mathbf{z}$.
Consequently, problem \eqref{eq:mlp_ref} has an objective function that is not jointly convex in $\mathbf{y}$ and $\mathbf{z}$, and the worst-case scenario may lie in the interior of the uncertainty set $\mathcal{Y}$.

Therefore, when assessing the validity of the feasibility time-step heuristic, it is necessary to consider not only the convexity of the operational problem but also whether the objective function \eqref{eq:mlp_obj} of the medial-level problem \eqref{medial_level} is jointly convex in the uncertain variables $\mathbf{y}$ and the operational variables $\mathbf{z}$ \citep{halemaneOptimalProcessDesign1983}.
The influence of the $3$ relevant nonconvexities in the operational problem and the joint nonconvexity of the objective function on robustness is summarized in Table \ref{tab:summary}.

In the following sections, we will describe the $3$ nonconvexities and identify the conditions under which they may pose challenges to the feasibility time-step heuristic.
\subsection{Piecewise-linear performance curves}
Nonlinear component efficiency curves can be approximated using piecewise-linear functions.
Rather than modeling the efficiency curves directly, it is often advantageous to describe the relationship between the energy inflow and outflow of a conversion component, as this relationship typically exhibits a more linear behavior \citep{vollAutomatedSuperstructurebasedSynthesis2013}.
The relevance of component performance curves in the context of robustness depends on whether the input energy carrier is externally supplied and unconstrained, or is subject to constraints or must be provided by the multi-energy system itself.
For example, meeting a cooling demand with a compression chiller induces an additional electricity demand within the system, whereas unconstrained externally supplied inputs, such as natural gas for a boiler, affect only operational costs.
Consequently, for robustness analysis, performance relations of components with unconstrained externally supplied inputs can be replaced by admissible energy outflow ranges.
In contrast, if a component induces additional demand for another energy carrier, the feasibility time-step heuristic may fail to provide a robust design.

The piecewise-linear relationship between the energy inflow and outflow of a conversion component $\Edot_{e_{in}}(\Edot_{e_{out}})$ can be modeled as follows according to \cite{sassModelCompendiumData2020}:
\begin{align*}
    & \Edot_{e_{in}} = \sum_{j \in \mathcal{J}} \left( b_{j} \lambda_{in, j} \frac{\Edot_{nom}}{\eta_{nom}} + \frac{\beta_{j}}{\eta_{nom}}\left(\Edot_{e_{out}, j} - \lambda_{out, j}  b_{j} \Edot_{nom}\right)\right)\\
    & \lambda_{out, j} \Edot_{nom} b_{j} \le \Edot_{e_{out}, j} & \forall j \in \mathcal{J} \\
    & \lambda_{out, j + 1} \Edot_{nom} b_{j} \ge \Edot_{e_{out}, j} & \forall j \in \mathcal{J} \\
    & \sum_{j \in \mathcal{J}} b_{j} = 1\\
    & \Edot_{e_{out}} = \sum_{j \in \mathcal{J}}{\Edot_{e_{out}, j}}
\end{align*}
Here, $\Edot_{e_{in}}$ and $\Edot_{e_{out}}$ denote the total energy inflow and outflow of the component, respectively, and $\Edot_{nom}$ is the nominal capacity of the component.
The binary variable $b_{j}$ selects the active segment $j \in \mathcal{J}$ of the piecewise-linear approximation.
The parameters $\lambda_{in, j}$ and $\lambda_{out, j}$ define the breakpoints of the approximation, and the slope
\begin{equation*}
    \beta_{j} = \frac{\lambda_{in, j + 1} - \lambda_{in, j}}{\lambda_{out, j + 1} - \lambda_{out, j}}
\end{equation*}
characterizes segment $j$.
The parameter $\eta_{nom}$ denotes the nominal efficiency.
This formulation represents a piecewise-linear approximation of the conversion efficiency.
During design optimization, the bilinear terms $b_{j} \Edot_{nom}$ can be linearized using Glover's reformulation \citep{glover1975}.
In the operational problem, $\Edot_{nom}$ is a fixed parameter, and the formulation directly reduces to a linear model.

A simple multi-energy system with two energy forms, $e_{in}$ and $e_{out}$, and allowed curtailment can be represented by the following operational feasibility problem:
\begin{align*}
    \underset{\mathbf{z} \in \mathcal{Z}, \phi \in \mathbb{R}}{\min} \quad & \phi\\
    s.t. \quad & \Edot_{e_{in}, agg} + \Edot_{e_{in}} \le \phi\\
    & \Edot_{e_{out}, agg} - \Edot_{e_{out}} \le \phi\\
    & \Edot_{e_{in}} = \sum_{j \in \mathcal{J}} b_{j} \lambda_{in, j} \frac{\Edot_{nom}}{\eta_{nom}} + \frac{\beta_{j}}{\eta_{nom}}\left(\Edot_{e_{out}, j} - \lambda_{out, j}  b_{j} \Edot_{nom}\right)\\
    & \lambda_{out, j} \Edot_{nom} b_{j} \le \Edot_{e_{out}, j} & \forall j \in \mathcal{J} \\
    & \lambda_{out, j + 1} \Edot_{nom} b_{j} \ge \Edot_{e_{out}, j} & \forall j \in \mathcal{J} \\
    & \sum_{j \in \mathcal{J}} b_{j} = 1\\
    & \Edot_{e_{out}} = \sum_{j \in \mathcal{J}}{\Edot_{e_{out}, j}}
\end{align*}
The component takes $\Edot_{e_{in}}$ as input energy and transforms it to $\Edot_{e_{out}}$, which it outputs, hence the opposite signs in the energy balance equations.
$\phi$ is an auxiliary variable introduced to model the energy gap, i.e.,
\begin{equation*}
    E_{gap}(\mathbf{x}, \mathbf{y}, \mathbf{z}) = \max \{\Edot_{e_{in}, agg} + \Edot_{e_{in}}, \Edot_{e_{out}, agg} - \Edot_{e_{out}}\},
\end{equation*}
and the operation is feasible if $\phi \le 0$.
For simplicity, all other components and demands have been aggregated into $\Edot_{e_{in}, agg}$ and $\Edot_{e_{out}, agg}$, respectively, to focus solely on the influence of the piecewise-linear energy inflow-outflow relationship.

Although many energy systems are modeled as MILP problems, the linearized performance curves are often convex and only appear nonconvex due to the modeling choice.
For example, the piecewise-linear functions used by \cite{vollAutomatedOptimizationbasedSynthesis2014} to model the energy inflow-outflow relation of a boiler, a combined heat and power unit, an absorption chiller, and a compression chiller are all convex functions.

If the linearized energy inflow-outflow curve is convex, a simplified formulation of $\Edot_{e_{in}}(\Edot_{e_{out}})$ can be used (cf. Chapter 1.3 in \cite{bertsimasIntroductionLinearOptimization1997} and Section 3.2.2 in \cite{varelmannSimultaneouslyOptimizingBidding2022}):
\begin{equation*}
    \Edot_{e_{in}} = \underset{j \in \mathcal{J}}{\max}\left\{b_{j} \lambda_{in, j} \frac{\Edot_{nom}}{\eta_{nom}} + \frac{\beta_{j}}{\eta_{nom}}\left(\Edot_{e_{out}} - \lambda_{out, j}  b_{j} \Edot_{nom}\right)\right\}
\end{equation*}
This results in the following problem formulation by dropping the constraints describing the piecewise-linear behavior in the general case:
\begin{align*}
    \underset{\mathbf{z} \in \mathcal{Z}, \phi \in \mathbb{R}}{\min} \quad & \phi\\
    s.t. \quad & \Edot_{e_{in}, agg} + \Edot_{e_{in}} \le \phi\\
    & \Edot_{e_{out}, agg} - \Edot_{e_{out}} \le \phi\\
    & \Edot_{e_{in}} = \underset{j \in \mathcal{J}}{\max}\left\{b_{j} \lambda_{in, j} \frac{\Edot_{nom}}{\eta_{nom}} + \frac{\beta_{j}}{\eta_{nom}}\left(\Edot_{e_{out}} - \lambda_{out, j}  b_{j} \Edot_{nom}\right)\right\}\\
    &\Edot_{e_{out}} \le \Edot_{nom}
\end{align*}
By substituting, we obtain the inequality
\begin{equation*}
    \Edot_{e_{in}, agg} +\underset{j \in \mathcal{J}}{\max}\left\{b_{j} \lambda_{in, j} \frac{\Edot_{nom}}{\eta_{nom}} + \frac{\beta_{j}}{\eta_{nom}}\left(\Edot_{e_{out}} - \lambda_{out, j}  b_{j} \Edot_{nom}\right)\right\} \le \phi,
\end{equation*}
which can be formulated as multiple linear constraints \citep{bertsimasIntroductionLinearOptimization1997, varelmannSimultaneouslyOptimizingBidding2022}:
\begin{equation}\label{eq:reformulation}
     \Edot_{e_{in}, agg} + b_{j} \lambda_{in, j} \frac{\Edot_{nom}}{\eta_{nom}} + \frac{\beta_{j}}{\eta_{nom}}\left(\Edot_{e_{out}} - \lambda_{out, j}  b_{j} \Edot_{nom}\right) \le \phi \quad \forall j \in \mathcal{J}
\end{equation}
Consequently, if the piecewise linearization of the energy inflow-outflow curve of the component is convex, i.e., $\Edot_{e_{in}}(\Edot_{e_{out}})$ is a convex function, it can be formulated fully linearly.
As a result, the designs identified by the feasibility time-step heuristic are guaranteed to be robust.

In contrast, if the energy inflow–outflow curve is nonconvex, this guarantee no longer holds; the feasibility time-step heuristic may still produce robust designs, but only incidentally.
For clarity, we relabel the energy forms by $1$ and $2$ instead of using $e_{in}$ and $e_{out}$ since input and output energies differ across components.
An example of a system where the feasibility time-step heuristic fails can be written as follows:
\begin{align*}
    \underset{\mathbf{x} \in \mathcal{X}, \mathbf{z} \in \mathcal{Z}}{\min} & \Edot_{max, c_{1}} + \Edot_{max, c_{2}}&&\\
    s.t. \quad & \Edot_{1, dem, d} + \Edot_{1, c_{1}, d} -  \Edot_{1, c_{2}, d} \le 0 \quad && \forall d \in \mathcal{D} \\
    & \Edot_{2, dem, d} - \Edot_{2, c_{1}, d} \le 0 \quad && \forall d \in \mathcal{D} \\
    & \Edot_{1, c_{1}, d} = \sum_{j \in \mathcal{J}} b_{c_{1}, d, j} \lambda_{in, j} \frac{\Edot_{nom}}{\eta_{nom}} + \frac{\beta_{j}}{\eta_{nom}}\left(\Edot_{2, c_{1}, d, j} - \lambda_{out, j}  b_{c_{1}, d, j} \Edot_{nom}\right) \quad && \forall d \in \mathcal{D} \\
    & \lambda_{out, j} \Edot_{nom} b_{c_{1}, d, j} \le \Edot_{2, c_{1}, d, j} \quad && \forall d \in \mathcal{D}, \forall j \in \mathcal{J} \\
    & \lambda_{out, j + 1} \Edot_{nom} b_{c_{1}, d, j} \ge \Edot_{2, c_{1}, d, j} \quad && \forall d \in \mathcal{D}, \forall j \in \mathcal{J} \\
    & \sum_{j \in \mathcal{J}} b_{c_{1}, d, j} = 1 \quad && \forall d \in \mathcal{D} \\
    & \Edot_{2, c_{1}, d} = \sum_{j \in \mathcal{J}}{\Edot_{2, c_{1}, d, j}} \quad && \forall d \in \mathcal{D} \\
    & \Edot_{1, c_{2}, d} \le \Edot_{max, c_{2}} \quad && \forall d \in \mathcal{D} \\
    & \Edot_{2, c_{1}, d} \le \Edot_{max, c_{1}} \quad && \forall d \in \mathcal{D} \\
\end{align*}
with
\begin{align*}
    &\mathbf{x} = [\Edot_{max, c_{1}}, \Edot_{max, c_{2}}] \\
    &\mathcal{X} = \mathbb{R}_{\ge 0}^{2} \\
    &\mathbf{z} = [\Edot_{1, c_{1}, 1}, \dots, \Edot_{1, c_{1}, |\mathcal{D}|}, \Edot_{2, c_{1}, 1}, \dots, \Edot_{2, c_{1}, |\mathcal{D}|}, \Edot_{1, c_{2}, 1}, \dots, \Edot_{1, c_{2}, |\mathcal{D}|}, \Edot_{2, c_{1}, 1, 1}, \dots, \Edot_{2, c_{1}, |\mathcal{D}|, |\mathcal{J}|},\\
    &b_{c_{1}, 1, 1}, \dots, b_{c_{1}, |\mathcal{D}|, |\mathcal{J}|}] \\
    &\mathcal{Z} = \mathbb{R}_{\ge 0}^{|\mathcal{D}| \left(3 + |\mathcal{J}|\right)} \times \{0, 1\}^{|\mathcal{D}| |\mathcal{J}|}
\end{align*}
\begin{figure}
    \centering
    \includegraphics[width=0.5\linewidth]{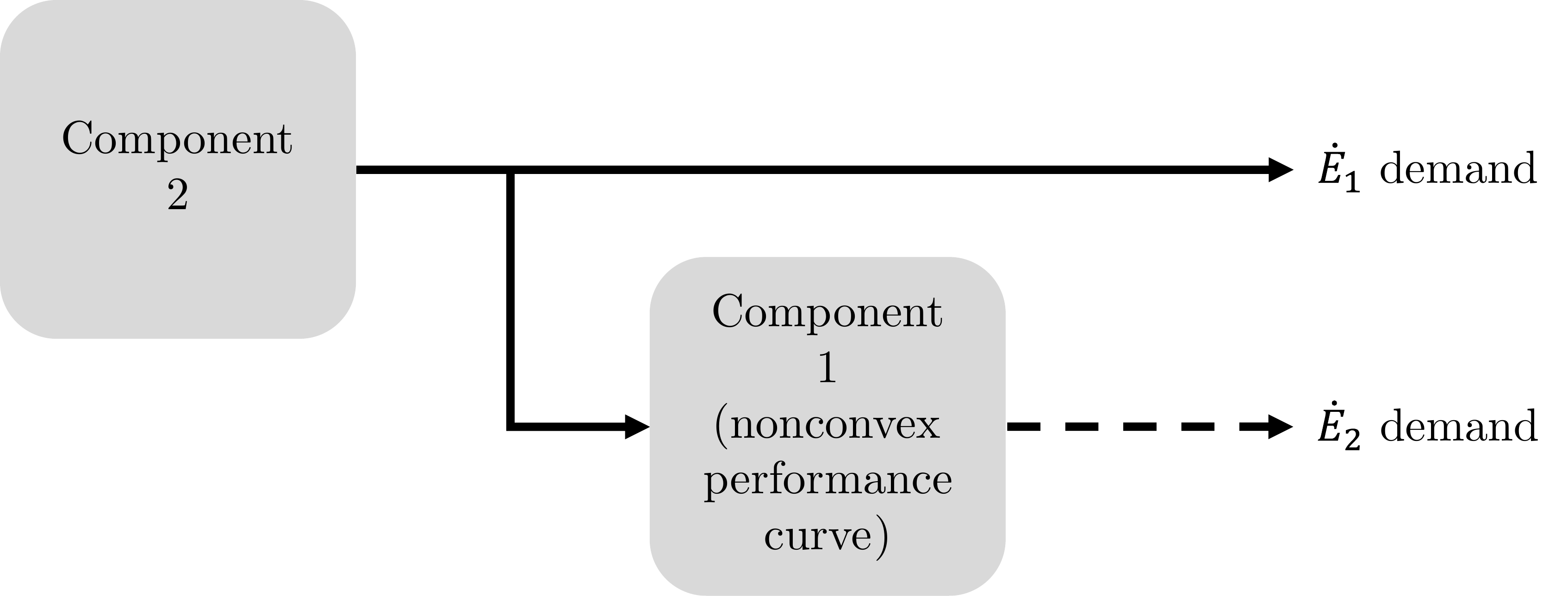}
    \caption{Illustration of an example system with a nonconvex piecewise-linear performance curve for component $1$. Component $2$ is supplied by an unconstrained external input energy, which is not shown here.}
    \label{fig:ex_non_convex}
\end{figure}
Here, $\mathbf{x}$ are the design variables and $\mathbf{z}$ are operational variables.
$\Edot_{max, c}$ are the component capacities, $c_{1}$ is a component with nonconvex $\Edot_{1}(\Edot_{2})$ and $c_{2}$ is an unconstrained externally supplied component that can supply $\Edot_{1}$.
Hence, the input energy of $c_{2}$ does not occur in any of the energy balances.
$\Edot_{1, c_{1}, d}$ and $\Edot_{2, c_{1}, d}$ denote the energy flows of energy forms $1$ and $2$ for component $1$ while $\Edot_{1, c_{2}, d}$ denotes the energy flow of energy form $1$ of component $2$ for data point $d$.
$b_{c_{1}, d, j}$ are the binary variables modeling the active segment of the piecewise linear energy inflow-outflow relationship at every data point $d$, and $\Edot_{2, c_{2}, d, j}$ are the energy flows of energy form $2$ for data point $d$ in segment $j$.
The parameters for the piecewise linear energy inflow-outflow relationship, namely $\lambda_{in, j}$, $\lambda_{out, j}$, $\beta_{j}$, and $\eta_{nom}$, as well as the demands $\Edot_{1, dem, d}$ and $\Edot_{2, dem, d}$ for energy forms $1$ and $2$ at each data point $d$, are provided in Section $1$ of the supplementary material.

\begin{figure}
    \centering
    \includegraphics[width=0.5\linewidth]{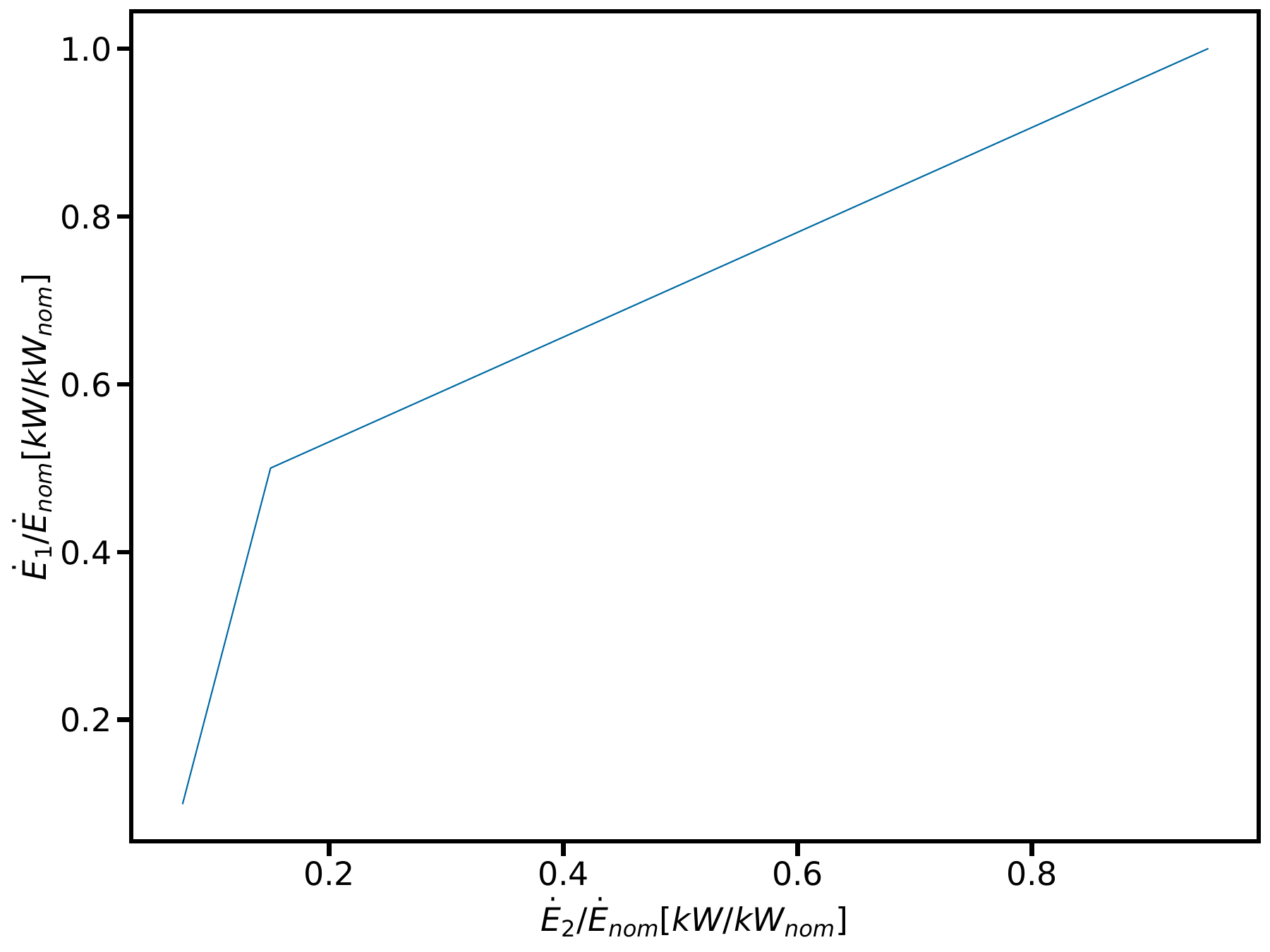}
    \caption{Nonconvex piecewise-linear performance curve for Component $c_{1}$ in Figure \ref{fig:ex_non_convex}.}
    \label{fig:ncvx_pwl}
\end{figure}

A schematic of the example system is shown in Figure \ref{fig:ex_non_convex} and the nonconvex piecewise-linear performance curve is illustrated in Figure \ref{fig:ncvx_pwl}.
The operational constraints are enforced for every data point $d$ in the historical data $\mathcal{D}$.

\begin{figure}
    \centering
    \includegraphics[width=0.5\linewidth]{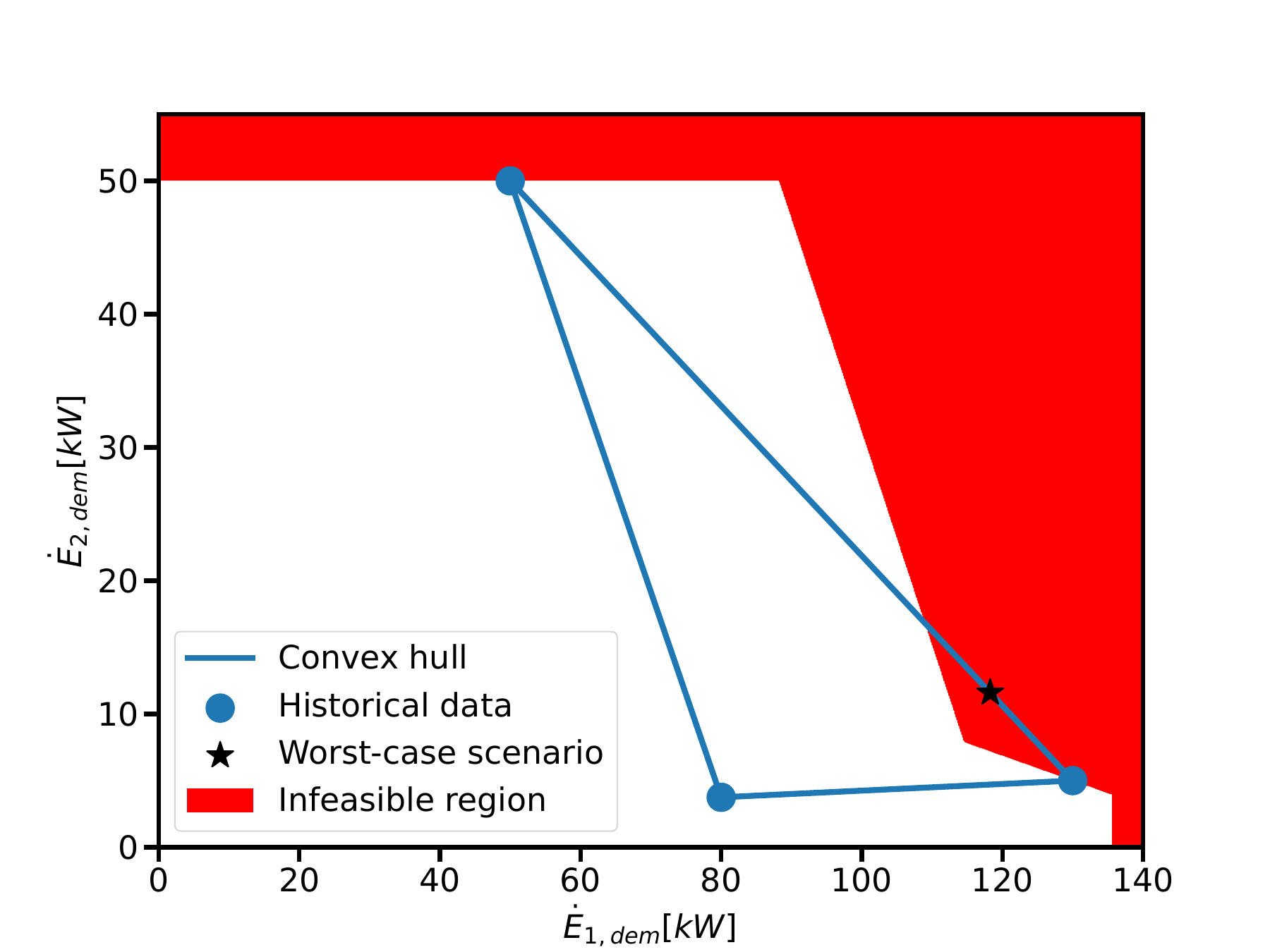}
    \caption{Counterexample that showcases how a design obtained using the feasibility time-step heuristic can fail to ensure feasibility within the convex hull of the historical data for a component with a nonconvex piecewise-linear energy inflow-outflow curve. The convex hull of historical demand data for energy forms $1$ and $2$ is plotted in blue. The infeasible region is shaded in red. The worst-case scenario inside the convex hull of the historical data is marked by a star. The example demonstrates that even if the heuristic identifies a design that is feasible for all historical scenarios, it may fail to guarantee feasibility for all points within the convex hull when nonconvexities are present.}
    \label{fig:pwl}
\end{figure}

Figure \ref{fig:pwl} illustrates the historical data and the feasible region of the resulting design.
Even though all three vertices of the convex hull are feasible, the convex hull contains infeasible regions.
The infeasible region is approximated by sampling uniformly spaced energy demands and solving the corresponding operational feasibility problem.
The worst-case scenario is identified by solving \eqref{medial_level} using an adaptive discretization algorithm \citep{blankenshipInfinitelyConstrainedOptimization1976}.

In the case where curtailment is not allowed, we have a semi-infinite equality constraint, which can be replaced by two corresponding inequality constraints, cf. Section \ref{sec:ncvx}, leading to the following \emph{operational} feasibility problem:
\begin{align*}
    \underset{\mathbf{z} \in \mathcal{Z}, \phi \in \mathbb{R}}{\min} \quad & \phi\\
    s.t. \quad & \Edot_{e_{in}, agg} + \Edot_{e_{in}} \le \phi\\
    & \Edot_{e_{out}, agg} - \Edot_{e_{out}} \le \phi\\
    & -\Edot_{e_{in}, agg} - \Edot_{e_{in}} \le \phi\\
    & - \Edot_{e_{out}, agg} + \Edot_{e_{out}} \le \phi\\
    & \Edot_{e_{in}} = \sum_{j \in \mathcal{J}} b_{j} \lambda_{in, j} \frac{\Edot_{nom}}{\eta_{nom}} + \frac{\beta_{j}}{\eta_{nom}}\left(\Edot_{e_{out}, j} - \lambda_{out, j}  b_{j} \Edot_{nom}\right)\\
    & \lambda_{out, j} \Edot_{nom} b_{j} \le \Edot_{e_{out}, j} & \forall j \in \mathcal{J} \\
    & \lambda_{out, j + 1} \Edot_{nom} b_{j} \ge \Edot_{e_{out}, j} & \forall j \in \mathcal{J} \\
    & \sum_{j \in \mathcal{J}} b_{j} = 1\\
    & \Edot_{e_{out}} = \sum_{j \in \mathcal{J}}{\Edot_{e_{out}, j}}\\
\end{align*}
Here, the inequalities with the opposite signs also appear:
\begin{equation*}
    E_{gap}(\mathbf{x}, \mathbf{y}, \mathbf{z}) = \max \{\Edot_{e_{in}, agg} + \Edot_{e_{in}}, \Edot_{e_{out}, agg} - \Edot_{e_{out}}, -\left(\Edot_{e_{in}, agg} + \Edot_{e_{in}}\right), -\left(\Edot_{e_{out}, agg} - \Edot_{e_{out}}\right)\}
\end{equation*}
Reformulation \eqref{eq:reformulation} used in the curtailment case is only applicable to constraints where the maximum term occurs with a positive sign (cf. Chapter 1.3 in \cite{bertsimasIntroductionLinearOptimization1997}).
$\Edot_{e_{in}}$, and hence the maximum term, occurs with a positive coefficient in one of the energy balances and with a negative coefficient in the complement of that energy balance, leading to the reformulation not being applicable in one case.
Consequently, the feasibility time-step heuristic is not guaranteed to yield robust designs, even if the piecewise-linear performance curves themselves are convex.

\subsection{Minimum part-load}
Many energy conversion units have a minimum part-load requirement, below which they cannot operate, and binary variables are commonly used to model the on/off behavior of such components \citep{mancoReviewMultiEnergy2024}:
\begin{align*}
    \Edot_{e_{out}} \le b \Edot_{max}\\
    b \Edot_{min} \le \Edot_{e_{out}}
\end{align*}
Here, $\Edot_{e_{out}}$ is the energy outflow of the component, and $\Edot_{min}$ and $\Edot_{max}$ are the minimum and maximum energy outflows, respectively.
If the binary variable $b = 0$, the component is off and $\Edot_{e_{out}} = 0$.
If $b = 1$, the component is on, and the minimum part-load is enforced, i.e., $\Edot_{min} \le \Edot_{e_{out}} \le \Edot_{max}$.
In our previous work \citep{wedemeyerRobustEnergySystem2025}, we demonstrated that when curtailment is not allowed, and the supplied energy outflow $\Edot_{e_{out}}$ must exactly meet the demand, minimum part-load constraints can render operation infeasible for low-demand scenarios.

\begin{figure}
    \centering
    \includegraphics[width=0.5\linewidth]{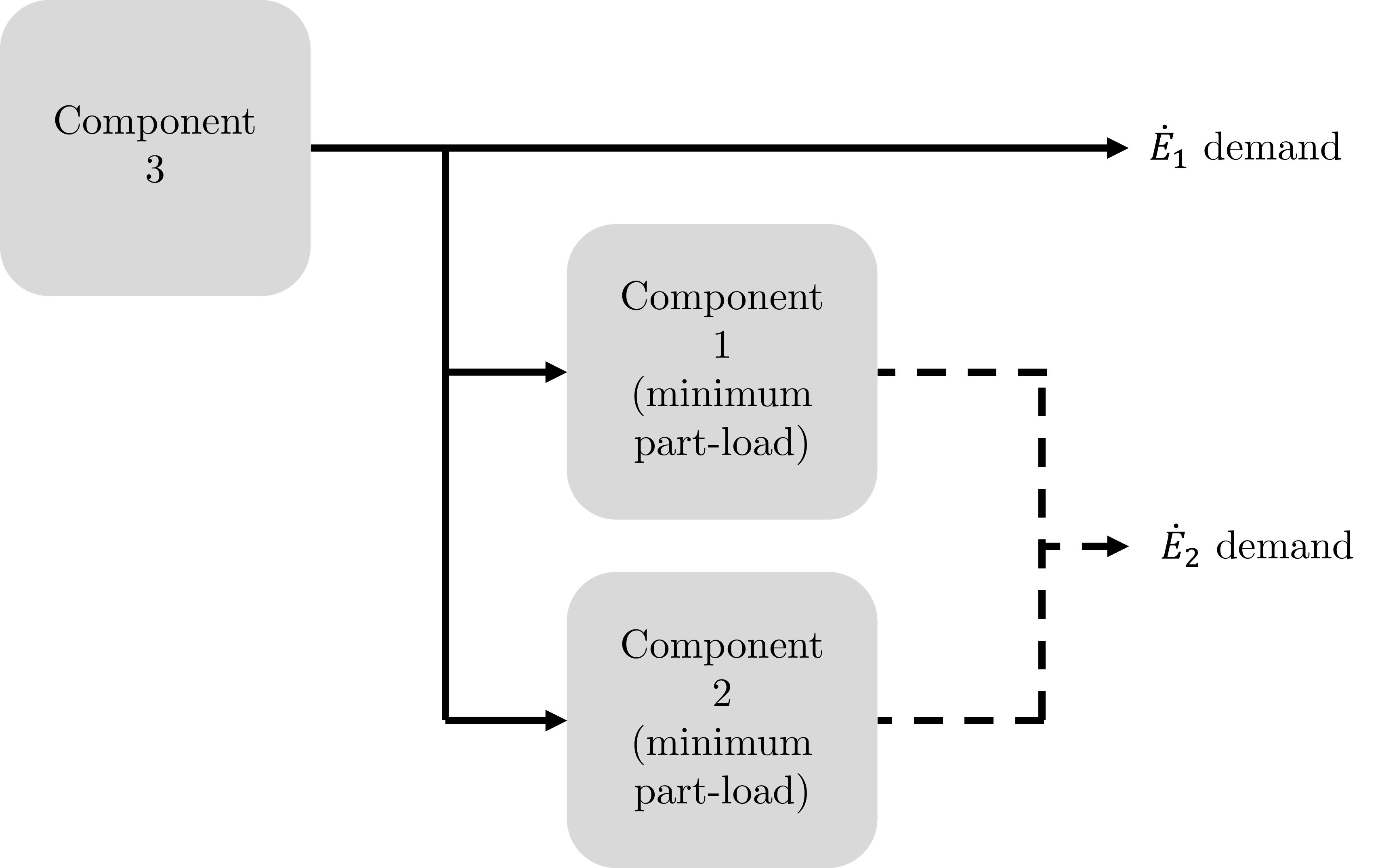}
    \caption{Illustration of the minimum part-load example system.}
    \label{fig:ex_min_part_load}
\end{figure}

However, if curtailment is allowed, the component can operate continuously at its minimum part-load, and any surplus energy can be dissipated.
Similar to the case of piecewise-linear performance curves, minimum part-load constraints can be replaced by admissible energy outflow ranges if the component has an unconstrained externally supplied energy inflow and hence cannot lead to the failure of the feasibility time-step heuristic.
In contrast, if the energy inflow induces an additional demand in a different energy form within the multi-energy system, operating at minimum part-load may create an extra demand that cannot be satisfied, potentially leading to infeasibility.
An example problem for which the feasibility time-step heuristic fails is shown in Figure \ref{fig:ex_min_part_load}.
It consists of three components.
Components $1$ and $2$ can convert energy form $1$ into energy form $2$, and component $3$ can supply energy form $1$ from an unconstrained externally supplied energy form, which is not shown in the figure.

The corresponding design problem can be formulated as
\begin{align*}
    \underset{\mathbf{x} \in \mathcal{X}, \mathbf{z} \in \mathcal{Z}}{\min} & \Edot_{max, c_{1}} + \Edot_{max, c_{2}} + \Edot_{max, c_{3}} &&\\
    \text{s.t.} \quad & \Edot_{1, dem, d} + \Edot_{1, c_{1}, d} + \Edot_{1, c_{2}, d} -  \Edot_{1, c_{3}, d} \le 0 \quad && \forall d \in \mathcal{D} \\
        & \Edot_{2, dem, d} - \Edot_{2, c_{1}, d} - \Edot_{2, c_{2}, d} \le 0 \quad && \forall d \in \mathcal{D} \\
        & \Edot_{1, c, d} = \Edot_{2, c, d} \quad && \forall c \in \{c_{1}, c_{2}\} \forall d \in \mathcal{D} \\
        & \Edot_{2, c, d} \le b_{c, d} \Edot_{max, c} \quad && \forall c \in \{c_{1}, c_{2}\} \forall d \in \mathcal{D}\\
        & b_{c, d} \Edot_{min, c} \le \Edot_{2, c, d} \quad && \forall c \in \{c_{1}, c_{2}\} \forall d \in \mathcal{D} \\
        & \Edot_{1, c_{3}, d} \le \Edot_{max, c_{3}} \quad && \forall d \in \mathcal{D} \\
        & \Edot_{min, c} = 0.2 \Edot_{max, c} \quad && \forall c \in \{c_{1}, c_{2}\} \\
        & \Edot_{max, c_{1}} \le 10 &&\\
        & 60 \le \Edot_{max, c_{2}} &&\\
\end{align*}
with
\begin{align*}
    &\mathbf{x} = [\Edot_{max, c_{1}}, \Edot_{max, c_{2}}, \Edot_{max, c_{3}}, \Edot_{min, c_{1}}, \Edot_{min, c_{2}}] \\
    &\mathcal{X} = \mathbb{R}_{\ge 0}^{5} \\
    &\mathbf{z} = [\Edot_{1, c_{1}, 1}, \dots, \Edot_{1, c_{1}, |\mathcal{D}|},  \Edot_{1, c_{2}, 1}, \dots, \Edot_{1, c_{2}, |\mathcal{D}|}, \Edot_{2, c_{1}, 1}, \dots, \Edot_{2, c_{1}, |\mathcal{D}|}, \Edot_{1, c_{3}, 1}, \dots, \Edot_{1, c_{3}, |\mathcal{D}|},\\
    &\Edot_{2, c_{2}, 1}, \dots, \Edot_{2, c_{2}, |\mathcal{D}|}, b_{c_{1}, 1}, \dots, b_{c_{1}, |\mathcal{D}|}, b_{c_{2}, 1}, \dots, b_{c_{2}, |\mathcal{D}|}] \\
    &\mathcal{Z} = \mathbb{R}_{\ge 0}^{5 |\mathcal{D}|} \times \{0, 1\}^{2 |\mathcal{D}|}
\end{align*}
where $\Edot_{max, c}$ are the component capacities.
The parameters for the demands can be found in Section $1$ of the supplementary material.
For simplicity, a constant conversion efficiency of $\eta = 1$ between energy inflow $\Edot_{1}$ and energy outflow $\Edot_{2}$ is assumed.
The operational constraints are enforced for every data point $d$ in the historical data $\mathcal{D}$.
Importantly, the allowed capacities for $c_{1}$ and $c_{2}$ differ and lead to a gap between the maximum output of component $1$ and the minimum part-load of component $2$.

\begin{figure}
    \centering
    \includegraphics[width=0.5\linewidth]{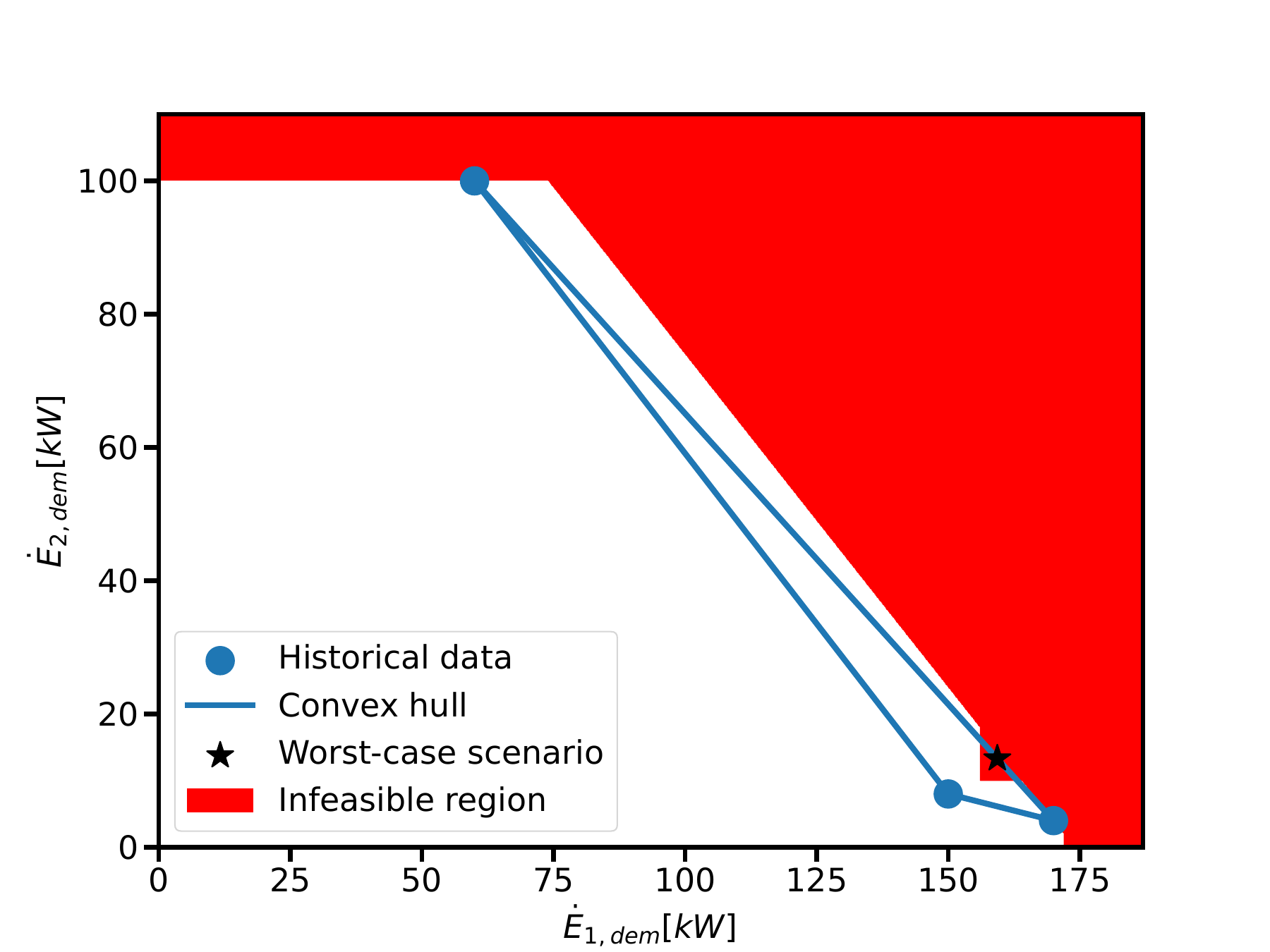}
    \caption{Counterexample for minimum part-load. The system schematic is shown in Figure \ref{fig:ex_min_part_load}. The figure illustrates how a design obtained using the feasibility time-step heuristic can be infeasible within the convex hull of the historical demand data. The convex hull of historical demand data for energy forms $1$ and $2$ is shown in blue, while the infeasible region is shaded in red. The worst-case scenario within the convex hull is marked by a star.}
    \label{fig:part_load}
\end{figure}

Figure \ref{fig:part_load} illustrates the convex hull of the historical dataset and the feasible region of the design identified by the introduced problem.
The infeasible region is approximated by sampling uniformly spaced energy demands and solving the corresponding operational feasibility problem.
The worst-case scenario is identified by solving \eqref{medial_level} using an adaptive discretization algorithm \citep{blankenshipInfinitelyConstrainedOptimization1976}.
Although all vertices of the convex hull are feasible, there are infeasible points inside the convex hull, introduced by the extra demand that is created when component $c_{2}$ has to be operated at minimum part-load.
\subsection{Storage complementarity}\label{sec:complementarity}
Storage complementarity refers to the restriction that an energy storage device cannot simultaneously charge and discharge, which, according to \cite{sassModelCompendiumData2020}, can be expressed as
\begin{align*}
    \Edot_{sto, in} \Edot_{sto, out} = 0,
\end{align*}
where $\Edot_{sto, in}$ is the storage energy inflow and $\Edot_{sto, out}$ is the storage energy outflow.
This nonconvex constraint is often applied when modeling battery systems \citep{appinoUseProbabilisticForecasts2018, sassModelCompendiumData2020}.
At the cell level, it is physically justified because the redox reaction proceeds in only one direction, depending on the applied potential \citep{goodenoughChallengesRechargeableLi2010}.
At the system level, simultaneous charging and discharging is also undesirable, as repeated cycling accelerates degradation \citep{kabirDegradationMechanismsLiion2017}.

Simultaneous charging and discharging can be interpreted as a form of curtailment, since the associated losses lead to energy dissipation.
In the following, we show that, if curtailment is allowed, storage complementarity does not affect the robustness of designs obtained using the feasibility time-step heuristic.
We achieve this by demonstrating that if storage complementarity is neglected during operational optimization, a feasible schedule with the same objective function that satisfies complementarity can be recovered in a post-processing step by adjusting the storage operation while preserving the trajectory of the storage level.
Since the operational problem is convex if storage complementarity is neglected, assuming no other nonconvexities are present, the design identified by the feasibility time-step heuristic is guaranteed to be robust.
Consequently, this design is then also robust for the same operational problem with enforced complementarity, as we can then construct a feasible complementarity-respecting schedule for every uncertainty realization.

In contrast, if curtailment is prohibited, this post-processing step is no longer applicable. Storage complementarity must then be enforced during optimization, introducing a nonconvexity that may result in non-robust designs.

We model the storage according to \cite{sassModelCompendiumData2020}.
The following constraints have to be satisfied by a feasible storage schedule:
\begin{align}
    & 0 \le \Edot_{sto, out, t} \le \Edot_{sto, out, max} \quad && \forall t \in \mathcal{T} \label{eq:out_bnds}\\
    & 0 \le \Edot_{sto, in, t} \le \Edot_{sto, in, max} \quad && \forall t \in \mathcal{T} \label{eq:in_bnds}\\
    & E_{sto, t}\left(1 + \frac{\Delta_{t}}{\tau}\right) - E_{sto, t - 1} = \Delta_{t} \left(\eta_{in} \Edot_{sto, in, t} - \frac{1}{\eta_{out}} \Edot_{sto, out, t}\right) \quad && \forall t \in \mathcal{T} \label{eq:sto} \\
    & \Edot_{agg, t} + \Edot_{sto, in, t} - \Edot_{sto, out, t}\le 0  \quad && \forall t \in \mathcal{T} \label{eq:balance}
\end{align}
Here, $E_{sto, t}$ denotes the storage level, $\Edot_{sto, in, t}$ and $\Edot_{sto, out, t}$ are the energy inflows and outflows of the storage, $t$ are the time steps and $\mathcal{T}$ is the set of all time steps, $\Delta_{t}$ is the time-step length, $\tau$ is the time constant for energy loss, and $\eta_{in}$ and $\eta_{out}$ are the charging and discharging efficiencies.
The differential equation describing the energy level of a generic storage is written as
\begin{equation*}
    \dfrac{\mathrm{d}E_{sto}}{\mathrm{d}t} = \eta_{in} \Edot_{sto, in} - \frac{1}{\eta_{out}} \Edot_{sto, out} - \frac{1}{\tau} E_{sto},
\end{equation*}
which we discretize using the implicit Euler scheme to obtain
\begin{equation*}
    E_{sto, t} =  E_{sto, t - 1} + \Delta_{t} \left(\eta_{in} \Edot_{sto, in, t} - \frac{1}{\eta_{out}} \Edot_{sto, out, t} - \frac{1}{\tau} E_{sto, t}\right)
\end{equation*}
and rearrange to obtain the equation describing the change of storage \eqref{eq:sto}.
The balance equation \eqref{eq:balance} ensures that the energy demand is met or exceeded.
$\Edot_{agg, t}$ aggregates the energy demand and the energy inflows and outflows of all other components and is introduced for ease of exposition.

Given a schedule $\left(\hat{E}_{sto, t}, \hat{\Edot}_{sto, in, t}, \hat{\Edot}_{sto, out, t}\right)$ that does not respect complementarity, a complementarity-respecting schedule can be constructed by determining $\left(\hat{E}_{sto, t}, \tilde{\Edot}_{sto, in, t}, \tilde{\Edot}_{sto, out, t}\right)$ such that the storage levels $\hat{E}_{sto, t}$ remain unchanged as is shown in the following.
Substituting into \eqref{eq:sto} yields the condition
\begin{equation}\label{eq:sto_ref}
    \Delta_{t}\left(\eta_{in} \hat{\Edot}_{sto, in, t} - \frac{1}{\eta_{out}} \hat{\Edot}_{sto, out, t} \right) = \Delta_{t}\left(\eta_{in} \tilde{\Edot}_{sto, in, t} - \frac{1}{\eta_{out}} \tilde{\Edot}_{sto, out, t} \right)
\end{equation}
Enforcing complementarity requires setting either $\tilde{\Edot}_{sto, in, t} = 0$ (discharging) or $\tilde{\Edot}_{sto, out, t} = 0$ (charging).
Which variable is set to zero is determined by evaluating the right-hand side of Equation \eqref{eq:sto}:
\begin{equation*}
    \eta_{in}\hat{\Edot}_{sto,in,t}
- \dfrac{1}{\eta_{out}}\hat{\Edot}_{sto,out,t}
\end{equation*}
Solving Equation \eqref{eq:sto_ref} for the remaining variable gives two cases:
\begin{equation*}
(\tilde{\Edot}_{sto,in,t}, \tilde{\Edot}_{sto,out,t})
=
\begin{cases}
\left(
0,\;
\hat{\Edot}_{sto,out,t}
- \eta_{in}\eta_{out}\hat{\Edot}_{sto,in,t}
\right),
& \text{if }
\eta_{in}\hat{\Edot}_{sto,in,t}
- \frac{1}{\eta_{out}}\hat{\Edot}_{sto,out,t}
\le 0, \\
\left(
\hat{\Edot}_{sto,in,t}
- \frac{1}{\eta_{in}\eta_{out}}\hat{\Edot}_{sto,out,t},\;
0
\right),
& \text{if }
\eta_{in}\hat{\Edot}_{sto,in,t}
- \frac{1}{\eta_{out}}\hat{\Edot}_{sto,out,t}
> 0
\end{cases}
\end{equation*}
First, consider the case
\begin{equation*}
    \eta_{in} \hat{\Edot}_{sto, in, t} - \frac{1}{\eta_{out}} \hat{\Edot}_{sto, out, t} \le 0
\end{equation*}
and
\begin{equation*}
   \tilde{\Edot}_{sto, out, t} = \hat{\Edot}_{sto, out, t} - \eta_{in} \eta_{out} \hat{\Edot}_{sto, in, t}
\end{equation*}
Feasibility is verified as follows. Since $\eta_{in} \eta_{out} \hat{\Edot}_{sto, in, t} \ge 0$,
\begin{equation*}
    \tilde{\Edot}_{sto, out, t} = \hat{\Edot}_{sto, out, t} - \underbrace{\eta_{in} \eta_{out} \hat{\Edot}_{sto, in, t}}_{\ge 0} \le \hat{\Edot}_{sto, out, t} \le \Edot_{sto, out, max}
\end{equation*}
Together with nonnegativity, i.e., $\tilde{\Edot}_{sto, out, t} = \left(\hat{\Edot}_{sto, out, t} - \eta_{in} \eta_{out} \hat{\Edot}_{sto, in, t} \right) \ge 0$, which follows from the case condition, Constraint \eqref{eq:out_bnds} is satisfied.
Trivially $\tilde{\Edot}_{sto, in, t} = 0$ satisfies Constraint \eqref{eq:in_bnds}.
It remains to verify the balance constraint \eqref{eq:balance}:
\begin{equation*}
    \Edot_{agg, t} - \tilde{\Edot}_{sto, out, t} = \Edot_{agg, t} - \left(\hat{\Edot}_{sto, out, t} - \eta_{in} \eta_{out} \hat{\Edot}_{sto, in, t} \right) \overset{!}{\le} \Edot_{agg, t} - \left(\hat{\Edot}_{sto, out, t} - \hat{\Edot}_{sto, in, t} \right) \le 0
\end{equation*}
To ensure $\Edot_{agg, t} - \left(\hat{\Edot}_{sto, out, t} - \eta_{in} \eta_{out} \hat{\Edot}_{sto, in, t} \right) \overset{!}{\le} \Edot_{agg, t} - \left(\hat{\Edot}_{sto, out, t} - \hat{\Edot}_{sto, in, t} \right)$, it is sufficient that $\eta_{in} \eta_{out} \le 1$.
For physically consistent storage systems, the charging and discharging efficiencies fulfill $0 \le \eta_{in} \le 1$ and $0 \le \eta_{out} \le 1$ and hence $\eta_{in} \eta_{out} \le 1$ holds.
Otherwise, the storage would generate energy when charging or discharging.
Therefore, the balance constraint \eqref{eq:balance} is satisfied.

Importantly, if curtailment is not allowed, the schedule needs to satisfy the balance constraint \eqref{eq:balance} with equality
\begin{equation*}
    \Edot_{agg, t} - \tilde{\Edot}_{sto, out, t} = \Edot_{agg, t} - \left(\hat{\Edot}_{sto, out, t} - \eta_{in} \eta_{out} \hat{\Edot}_{sto, in, t} \right) \overset{!}{=} \Edot_{agg, t} - \left(\hat{\Edot}_{sto, out, t} - \hat{\Edot}_{sto, in, t} \right) = 0,
\end{equation*}
which is only the case for $\eta_{in} \eta_{out} = 1$.
Hence, for $\eta_{in} < 1$ or $\eta_{out} < 1$, complementarity must be enforced in the operational problem, which introduces a nonconvexity that may lead to designs identified by the feasibility time-step heuristic not being robust.

Now consider the case
\begin{equation*}
    \eta_{in} \hat{\Edot}_{sto, in, t} - \frac{1}{\eta_{out}} \hat{\Edot}_{sto, out, t} > 0
\end{equation*} and
\begin{equation*}
    \tilde{\Edot}_{sto, in, t} = \hat{\Edot}_{sto,in,t} - \frac{1}{\eta_{in}\eta_{out}}\hat{\Edot}_{sto,out,t}
\end{equation*}
Feasibility follows similarly:
\begin{equation*}
    \tilde{\Edot}_{sto, in, t} = \hat{\Edot}_{sto, in, t} - \underbrace{\frac{1}{\eta_{in} \eta_{out} }\hat{\Edot}_{sto, out, t}}_{\ge 0} \le \hat{\Edot}_{sto, in, t} \le \Edot_{sto, in, max}
\end{equation*}
Together with nonnegativity, which follows from the case condition, Constraint \eqref{eq:in_bnds} is satisfied.
Trivially $\tilde{\Edot}_{sto, out, t} = 0$ satisfies Constraint \eqref{eq:out_bnds}.
Again, the balance constraint must hold:
\begin{equation*}
    \Edot_{agg, t} + \tilde{\Edot}_{sto, in, t} = \Edot_{agg, t} + \left(\hat{\Edot}_{sto, in, t} - \frac{1}{\eta_{in} \eta_{out} }\hat{\Edot}_{sto, out, t}\right) \overset{!}{\le} \Edot_{agg, t} + \left(\hat{\Edot}_{sto, in, t} - \hat{\Edot}_{sto, out, t} \right) \le 0,
\end{equation*}
which again holds true if $\eta_{in} \eta_{out} \le 1$.

In summary, as long as $\eta_{in} \eta_{out} \le 1$, which is the physically relevant case, storage complementarity can be neglected if curtailment is allowed, as a feasible complementarity-respecting schedule can always be constructed in post-processing, guaranteeing that designs identified by the feasibility time-step heuristic are robust.
\subsection{Guidelines and mitigating factors}
Table \ref{tab:summary} summarizes our results from this section and indicates the cases in which the feasibility time-step heuristic is guaranteed to identify robust designs.
\begin{table}[h]
    \begin{center}
    \begin{tabular}{l|c|c}
        \hline
        Nonconvexity type & Curtailment & No-curtailment \\
        \hline
        Piecewise-linear convex energy inflow-outflow curve& \cmark & \xmark \\
        Piecewise-linear nonconvex energy inflow-outflow curve& \xmark & \xmark \\
        Minimum part-load & \xmark & \xmark \\
        Storage complementarity & \cmark & \xmark \\
        Objective function not jointly convex in $\mathbf{y}$ and $\mathbf{z}$ & \xmark & \xmark \\
        \hline
    \end{tabular}
    \end{center}
    \caption{Overview of the investigated types of nonconvexities: A checkmark indicates that the feasibility time-step heuristic is guaranteed to yield a robust design despite the nonconvexity, whereas a cross indicates that the design may not be robust.}
    \label{tab:summary}
\end{table}

If the objective function of \eqref{medial_level} is not jointly convex in the uncertainties and the operational variables, the feasibility time-step heuristic does not, in general, guarantee robust designs.
In the absence of curtailment, any nonconvexity can lead to solutions that are not robust.
In contrast, if curtailment is allowed, storage complementarity and piecewise-linear convex energy inflow-outflow curves do not compromise the robustness of designs identified by the heuristic.
However, nonconvex energy inflow-outflow curves or minimum part-loads may lead the feasibility time-step heuristic to fail.
We have shown small single-time-step examples for which the design identified by the feasibility time-step heuristic is infeasible within the convex hull of historical data.

The theoretical analysis in this section relies on the assumption that unlimited curtailment is allowed.
In practical settings, however, curtailment may be limited by an upper bound.
If curtailment is allowed but bounded, the robustness guarantees derived for the unlimited-curtailment case no longer apply.
The definition of the energy gap then becomes:
\begin{equation*}
    \begin{aligned}
    &E_{gap} (\mathbf{x}, \mathbf{y}, \mathbf{z}) = &\\
    &\underset{ e \in \mathcal{E}, t \in \mathcal{T}}{\max} \biggl\{ \max && \Bigl\{ \Edot_{e, dem, t} - \sum_{c \in \mathcal{C}_{e, out}} \Edot_{e, c, t} + \sum_{c \in \mathcal{C}_{e, in}} \Edot_{e, c, t} + \sum_{c \in \mathcal{C}_{e, sto}} \left(\Edot_{e, c, t, in} - \Edot_{e, c, t, out} \right),\\
    &&&-\Edot_{e, dem, t} + \sum_{c \in \mathcal{C}_{e, out}} \Edot_{e, c, t} - \sum_{c \in \mathcal{C}_{e, in}} \Edot_{e, c, t} - \sum_{c \in \mathcal{C}_{e, sto}} \left(\Edot_{e, c, t, in} - \Edot_{e, c, t, out} \right) - \Edot_{e, t, cur} \Bigr\} \biggr\},
    \end{aligned}
\end{equation*}
where $\Edot_{e, t, cur}$ is the maximum curtailment of energy form $e$ at time $t$.
As in the no-curtailment case, the energy flow terms enter with opposite signs in the two expressions inside the maximum.
Therefore, the reformulations used for the unlimited-curtailment case are no longer directly applicable.
However, due to the additional curtailment term, the resulting problem can be formulated as a regular ESIP \citep{djelassiGlobalSolutionSemiinfinite2021}.
Hence, unlike the problems arising in the no-curtailment case, it can be solved using adaptive discretization algorithms \citep{djelassidiskretisierungsbasiertealgorithmenfur2020}.

If storage components are included in the design, coupling between time steps increases system flexibility by allowing energy to be shifted across time steps.
In the minimum part-load case, storage may compensate for low-demand periods where the installed components cannot operate at or below their minimum part-load by shifting generated energy flows from other time steps.
Similarly, energy losses in storage by simultaneous charging and discharging (e.g., thermal energy storage) can act as a weak form of curtailment.
However, this curtailment is limited and hard to estimate a priori.

In general, the feasibility time-step heuristic can always be used to identify candidate designs.
Nevertheless, when the problematic nonconvexities identified in Table \ref{tab:summary} are present, verifying the feasibility of the design is recommended.
At present, this can only be achieved by solving \eqref{medial_level} to global optimality given the identified design, an approach that may quickly become computationally intractable.
\section{Illustrative energy system example}\label{sec:case_study}
\begin{figure}
    \centering
    \includegraphics[width=0.9\linewidth]{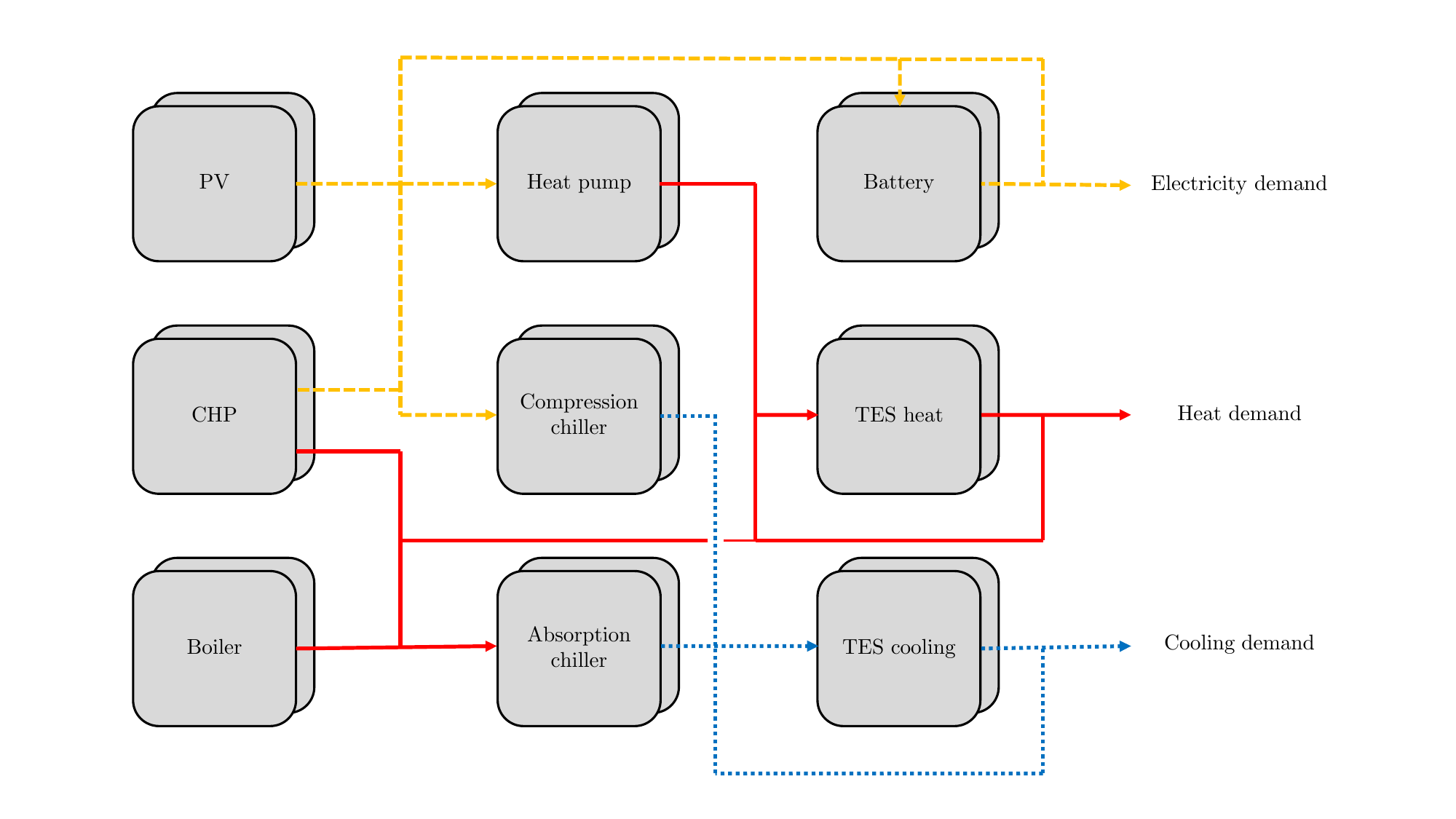}
    \caption{Illustration of the modified example system based on \cite{sassModelCompendiumData2020}. Here, no connection to the electrical grid is present. Up to two copies of each component may be installed. For the CHPs, three categories are available: small, medium, and large, with different nominal capacity ranges and economic data. Note that the gas inflow into the CHP and boiler components is omitted from the illustration, as gas is externally supplied and unconstrained. Thermal energy storage units are abbreviated by TES.}\label{fig:sass_case_study}
\end{figure}

We investigate an illustrative example based on the work of \cite{sassModelCompendiumData2020} and \cite{vollAutomatedSuperstructurebasedSynthesis2013}, which contains most of the previously discussed nonconvexities.
The only nonconvexity not present is a nonconvex energy inflow-outflow curve, as we were unable to identify components with such curves in the literature.
We modify the case study from \cite{sassModelCompendiumData2020} by removing the connection to the electricity grid, leading to an island multi-energy system; the only imported energy source is gas, used to power the CHPs and boilers.
Figure \ref{fig:sass_case_study} showcases the investigated system.
Up to two copies of each component may be installed.
For the CHPs, three categories are available: small, medium, and large, with different nominal capacity ranges and economic data.
Unlike \cite{sassModelCompendiumData2020}, we only investigate the minimization of the total annualized cost (TAC) of the energy system and do not consider the global warming impact, as our focus lies on investigating the robustness of the identified designs.

The data for the demands, solar irradiation, and ambient temperature are obtained from \cite{sassModelCompendiumData2020} and accessible at \cite{HECIEnergybenchmarkGitLab2020}.
To construct representative days from the historical data, we follow the data preprocessing methodology of our previous work \citep{wedemeyerRobustEnergySystem2025}.
The data are separated into days and resampled to $12$ time steps per day, then normalized to have unit variance and zero mean, and clustered into $4$ representative scenarios using k-means clustering \citep{macqueen1967some}.
All components are modeled according to \cite{sassModelCompendiumData2020}.
The investigated system contains components with piecewise-linear performance curves.
We observed that the piecewise-linear performance curves provided by \cite{sassModelCompendiumData2020} are nonconvex and are inconsistent with the provided original nonlinear curves.
Furthermore, they also disagree with the piecewise linearizations provided by \cite{vollAutomatedSuperstructurebasedSynthesis2013} on which the work of \cite{sassModelCompendiumData2020} is based.
Hence, we re-linearized the nonlinear performance curves provided by \cite{sassModelCompendiumData2020}.
For our linearization method and results, see Section 2 in the supplementary material.
Additionally, the linearization data from \cite{sassModelCompendiumData2020} suggests different minimum part-loads for the electricity and heat output of the CHPs.
We follow the approach of \cite{vollAutomatedSuperstructurebasedSynthesis2013}, who assume a shared minimum part-load of $0.5$.

For a complete description of the parameters used and the optimization model, we refer to Section 3 of the supplementary material.

The linearized energy inflow-outflow curves of the boiler, CHP, and heat pump only have a single linear segment and are hence convex.
The PV system and the energy storage systems are modeled with constant efficiencies.
Only the linearized energy inflow-outflow curves of compression and absorption chillers are modeled with two linear segments each, and the curves are also convex.
The battery is modeled with charging complementarity constraints.
However, those can be neglected with respect to robustness as described in Section \ref{sec:complementarity}.
Regarding joint convexity of the objective function, the nominal efficiency of the heat pump depends on the uncertain ambient temperature and enters the lower-level problem as follows:
\begin{equation*}
    \Edot_{e_{in}, hp, t} - \frac{\Edot_{e_{out}, hp, t}}{\eta_{hp, nom, t}(T_{amb, t})} = 0
\end{equation*}
To remove the coupling equality constraint, we directly insert this definition of the energy inflow into the energy balance, cf. Section \ref{sec:ncvx}.
The resulting term is not jointly convex in $\mathbf{y}$ and $\mathbf{z}$ and hence may lead to a failure of the feasibility time-step heuristic \citep{teichgraeberExtremeEventsTime2020}.

Consequently, from the preceding analysis, the only sources of nonconvexity in this case study that may cause the feasibility time-step heuristic to fail are the objective function not being jointly convex and the presence of minimum part-load constraints.
We determine a candidate robust design using the feasibility time-step heuristic \citep{teichgraeberExtremeEventsTime2020} and investigate its robustness.
\subsection{Rigorously verifying robustness}\label{sec:method}
To validate the robustness of the candidate design generated by the feasibility time-step heuristic \citep{teichgraeberExtremeEventsTime2020} and to determine whether minimum part-load constraints or joint nonconvexities in the objective function introduce infeasibilities, we solve problem \eqref{medial_level} to rigorously identify the maximum violation within the convex hull of the historical data.

\eqref{medial_level} is a challenging bilevel problem, and preliminary numerical experiments showed prohibitively long convergence times when using adaptive discretization approaches \citep{blankenshipInfinitelyConstrainedOptimization1976}.
Furthermore, the presence of integer variables renders the lower-level problem nonconvex, making common single-level reformulation approaches based on the KKT-optimality conditions \citep{fortuny-amatRepresentationEconomicInterpretation1981} or the linear dual \citep{burgardOptknockBilevelProgramming2003} inapplicable.
To address this issue, we propose a hybrid approach that combines elements of classical adaptive discretization methods \citep{blankenshipInfinitelyConstrainedOptimization1976} with a single-level reformulation by embedding the linear dual of the lower-level problem \citep{burgardOptknockBilevelProgramming2003}.
Specifically, we apply an adaptive discretization approach to the integer variables and then embed the dual of the resulting linear program into the medial-level problem for each set of discretized integer variables.

Our approach is related to the approach proposed by \cite{zhaoExactAlgorithmTwostage}.
They consider a similar decomposition of a MILP lower-level problem: The continuous lower-level problem is dualized for fixed integer variables, while the integer variables are handled through column-and-constraint generation \citep{zengSolvingTwostageRobust2013}.
A similar approach has been proposed by \cite{mitsosGlobalSolutionNonlinear2010} to solve a mixed-integer bilevel program, where the KKT conditions are used to tighten the lower bounding scheme.
In contrast, we do not introduce a parametric upper bound on the objective value of the lower-level problem and instead embed the linear dual at each discretization point.
\begin{align*}\tag{LLP}\label{lower_level}
     \underset{\mathbf{z}_{c} \in \mathbb{R}^{n_{c}}, \mathbf{z}_{d} \in \{0, 1\}^{n_{d}}}{\min} & \quad  \mathbf{c}_{c}^T \mathbf{z}_{c} \\
     \text{s.t.} & \quad \mathbf{A}_{c}(\mathbf{y}) \mathbf{z}_{c} \le \mathbf{b}(\mathbf{y}) - \mathbf{A}_{d} \mathbf{z}_{d}
\end{align*}
To this end, we first introduce how we form the dual of the lower-level problem with fixed integer variables that is then embedded into the medial-level problem.
\eqref{lower_level} shows the lower-level problem, rewritten in inequality-constrained form without explicit bounds on the continuous variables, which have been moved to the constraints, to simplify the dual formulation.
For a general treatment of dual problems of linear programs, see \cite{bertsimasIntroductionLinearOptimization1997}.
Here, $\mathbf{c}_{c}$ is the cost vector and $\mathbf{A}_{c}(\mathbf{y})$ and $\mathbf{A}_{d}$ are the coefficient matrices for the linear constraints for the continuous and integer variables, respectively.
$\mathbf{A}_{c}(\mathbf{y})$ may be a function of the medial-level variables.
For example, in the present case study, the uncertain variables include ambient temperatures $T_{amb, t}$, which affect the efficiency of the heat pumps: $\eta_{hp, nom, t}(T_{amb, t}) = \frac{0.36T_{hp}}{T_{hp} - T_{amb, t}}$, where $T_{hp}$ is the output temperature of the heat pump \citep{sassModelCompendiumData2020} and may appear in the lower-level problem in constraints as $\Edot_{e_{in}, hp, t} - \frac{\Edot_{e_{out}, hp, t}}{\eta_{hp, nom, t}} = 0$.
The vector $\mathbf{b}(\mathbf{y})$ denotes the right-hand side of the constraints and depends linearly on the medial-level variables $\mathbf{y}$.

Importantly, we assume here that \eqref{lower_level} corresponds to the lower level of an SIP or has been transformed into an SIP by moving the coupling constraints into the objective function, as described in Section \ref{sec:PS}.
Under this assumption, the dependence of $\mathbf{A}_{c}(\mathbf{y})$ and $\mathbf{b}(\mathbf{y})$ on the medial-level variables $\mathbf{y}$ results solely from the reformulation of the nonlinear maximum terms in the objective function of the medial-level problem \eqref{eq:mlp_obj} into an equivalent linear form (cf. Chapter 1.3 in \cite{bertsimasIntroductionLinearOptimization1997} and Section 3.2.2 in \cite{varelmannSimultaneouslyOptimizingBidding2022}).
Hence, the feasibility of the lower-level problem does not depend on $\mathbf{y}$.

When forming the dual, the integer variables $\mathbf{z}_{d}$ are fixed and therefore appear as constants on the right-hand side.
Forming the dual of problem \eqref{lower_level} with fixed integer variables $\bar{\mathbf{z}}_{d}$ yields problem \eqref{lower_level_dual}.
\begin{align*}\tag{LLP dual}\label{lower_level_dual}
     \underset{\boldsymbol{\lambda} \in \mathbb{R}_{\le 0}^{n_{cs}}}{\max} & \quad \boldsymbol{\lambda}^T \left(\mathbf{b}(\mathbf{y}) - \mathbf{A}_{d} \bar{\mathbf{z}}_{d}\right) \\
     \text{s.t.} & \quad \boldsymbol{\lambda}^T  \mathbf{A}_{c}(\mathbf{y}) = \mathbf{c}_{c}^{T}
\end{align*}

In contrast to classical adaptive discretization methods, which discretize all lower-level variables \citep{blankenshipInfinitelyConstrainedOptimization1976, djelassiGlobalSolutionSemiinfinite2021}, we discretize only the integer variables $\mathbf{z}_{d}$.
Let $\bar{\mathbf{z}}_{d, k}$ denote a fixed realization of the integer variables corresponding to the discretization index $k \in \mathcal{K}$ from the discretization set $\mathcal{Z}_{disc, d}$.
The resulting discretized medial-level problem can be written as
\begin{align*}\tag{MLP-disc}\label{medial_level_disc}
     \underset{\mathbf{y} \in \mathcal{Y}, \phi \in \mathbb{R}}{\max} & \quad \phi \\
     \text{s.t.} & \quad \phi \le \underset{\mathbf{z}_{c, k} \in \mathcal{Z}_{c}(\mathbf{y}, \bar{\mathbf{z}}_{d, k})}{\min} \quad \mathbf{c}_{c}^T \mathbf{z}_{c, k} \ \forall k \in \mathcal{K}
\end{align*}

Here, $\phi$ is an auxiliary variable that ensures that the objective is smaller than the optimal value of any of the discretizations so far, and 
$\mathcal{Z}_{c}(\mathbf{y}, \bar{\mathbf{z}}_{d, k}) = \left\{\mathbf{z}_{c} \ | \  \mathbf{A}_{c}(\mathbf{y}) \mathbf{z}_{c} \le \mathbf{b}(\mathbf{y}) - \mathbf{A}_{d} \bar{\mathbf{z}}_{d, k}, \mathbf{z}_{c} \in \mathbb{R}^{n_{c}} \right\}$ denotes the feasible set for the continuous lower-level variables.

Since the objective here is to verify the robustness of a given design, the design variables $\mathbf{x}$ are treated as fixed and their dependence is omitted from the notation for simplicity, unlike in problem \eqref{medial_level}.

\eqref{medial_level_disc} is a bilevel problem.
To obtain a single-level formulation, we replace the embedded lower-level problem by its dual for each discretization point:
\begin{align*}\tag{MLP-dual-disc}\label{medial_level_dual_disc}
     \underset{\mathbf{y} \in \mathcal{Y}, \phi \in \mathbb{R}, \boldsymbol{\lambda}_{k} \in \mathbb{R}_{\le 0}^{n_{cs}} \forall k \in \mathcal{K}, \mathbf{z}_{c, k} \in \mathbb{R}^{n_{c}} \forall k \in \mathcal{K}}{\max} & \quad \phi &\\
     \text{s.t.} \quad & \phi \le \boldsymbol{\lambda}_{k}^T \left(\mathbf{b}(\mathbf{y}) - \mathbf{A}_{d} \bar{\mathbf{z}}_{d, k}\right) &\ \forall k \in \mathcal{K} \\
     & \boldsymbol{\lambda}_{k}^T  \mathbf{A}_{c}(\mathbf{y}) = \mathbf{c}_{c}^{T} & \ \forall k \in \mathcal{K}\\
     & \mathbf{A}_{c}(\mathbf{y}) \mathbf{z}_{c, k} \le \mathbf{b}(\mathbf{y}) - \mathbf{A}_{d} \bar{\mathbf{z}}_{d, k} &\ \forall k \in \mathcal{K}\\
     & \boldsymbol{\lambda}_{k}^T \left(\mathbf{b}(\mathbf{y}) - \mathbf{A}_{d} \bar{\mathbf{z}}_{d, k}\right) = \mathbf{c}_{c}^T \mathbf{z}_{c, k} &\ \forall k \in \mathcal{K}
\end{align*}
Note that we additionally embed the primal constraints and set the primal and dual objective equal, as we found this to improve performance in preliminary studies.

Problem \eqref{medial_level} is solved by adaptive discretization of the integer lower-level variables, following the idea of \cite{blankenshipInfinitelyConstrainedOptimization1976}.
In each iteration, a worst-case candidate $\mathbf{y}$ is determined by solving problem \eqref{medial_level_dual_disc}.
Since \eqref{medial_level_disc}, and hence \eqref{medial_level_dual_disc}, is a relaxation of \eqref{medial_level}, it provides an upper bound on the objective value of \eqref{medial_level}.
\eqref{medial_level_dual_disc} is a relaxation of \eqref{medial_level} because the constraint bounding the auxiliary variable $\phi$ in \eqref{medial_level_dual_disc} is only enforced for the discretization points, which are a subset of all possible integer variable realizations.

In the following, we assume that \eqref{lower_level} is feasible for every $\mathbf{y} \in \mathcal{Y}$. This is satisfied in the considered problem class because the feasibility of the lower-level problem does not depend on $\mathbf{y}$. We further assume that \eqref{lower_level} is bounded. Since its objective is the energy gap, an unbounded lower-level problem would imply that at least one energy flow can become arbitrarily large, which would indicate a missing physical or operational bound and hence a modeling error.

The lower-level problem \eqref{lower_level} is subsequently solved for a fixed worst-case candidate $\mathbf{y}$ to obtain a new realization of the integer variables, which is added to the discretization set.
This solution provides a feasible point for \eqref{medial_level} and therefore a lower bound on its optimal objective value. The worst-case candidate with the highest objective value, i.e., the highest lower bound, is retained as the incumbent for the globally optimal worst-case scenario.

We initialize the discretization set with the integer variable values of the solution to the operational problem for the historical data point with the highest objective value, but additional initial discretization points may also be used.
The steps described above are then repeated.
Whenever solving \eqref{medial_level_dual_disc} or \eqref{lower_level} yields an improvement to the current upper bound $f_{ub}$ or lower bound $f_{lb}$, respectively, the corresponding bound is updated.
The procedure terminates when the upper and lower bounds differ by at most the specified optimality tolerance $\epsilon_{o}$.
Consequently, the identified worst-case scenario is certified to be globally $\epsilon_{o}$-optimal.
Additionally, if the upper bound is less than or equal to a user-specified feasibility tolerance $\epsilon_{f}$, the worst-case scenario is deemed to be feasible since the largest value of the energy gap \eqref{eq:mlp_obj} is at most the feasibility tolerance.
In this case, the incumbent worst-case scenario is not globally $\epsilon_o$-optimal.
Depending on the problem formulation, the energy forms included in the definition of the energy gap may differ in scale.
Consequently, the resulting worst-case scenario may be dominated by a particular energy form.
The energy gap is defined in terms of the largest absolute deviation because its primary purpose is to detect whether any infeasibility occurs, rather than to quantify or compare the relative impact of individual infeasibilities.
In particular, a negative energy gap guarantees system feasibility regardless of the relative magnitudes of the individual energy contributions.

To balance the influence of the different energy forms, their contributions could be normalized or multiplied by suitable positive scaling factors.
Such a rescaling would retain the property that a negative energy gap implies feasibility.
In the present work, however, we do not apply any normalization or scaling.

Let $k_{0}$ denote the number of initial discretization points, and let
$\bar{\mathbf{z}}_{d,k}$, $k = 1,\dots,k_{0}$, be the corresponding fixed initial realizations of the integer lower-level variables.
The algorithm to rigorously verify the robustness of a candidate design then reads:
\begin{algorithm}[H]
\caption{B\&F-based verification algorithm}
\label{alg:ver}
\begin{algorithmic}
\State $f_{lb} \gets -\infty$, $f_{ub} \gets \infty$, $\mathcal{K} \gets \{1, \dots, k_{0}\}$, $k_{last} \gets k_{0}$
\While{$(f_{ub} - f_{lb}) \ge \epsilon_{o} \wedge f_{ub} \ge \epsilon_{f}$}
\State solve \eqref{medial_level_dual_disc}
\State $f_{ub} \gets \phi^{*}$
\State solve \eqref{lower_level}
\If{$\mathbf{c}_{c}^T \mathbf{z}_{c}^{*} \geq f_{lb}$}
\State $f_{lb} \gets \mathbf{c}_{c}^T \mathbf{z}_{c}^{*}$
\EndIf
\State $k_{last} \gets 1 + k_{last}$
\State $\bar{\mathbf{z}}_{d, k_{last}} \gets \mathbf{z}_{d}^{*}$
\State $\mathcal{K} \gets \mathcal{K} \cup \{k_{last}\}$
\EndWhile
\end{algorithmic}
\end{algorithm}

\cref{alg:ver} is the algorithm introduced by \cite{blankenshipInfinitelyConstrainedOptimization1976} with the exception that the continuous variables of the lower-level problem have been lifted into the upper-level problem.
As described previously, for our problem class, the lower-level problem is always feasible and bounded, hence there always exists an optimal solution, and the embedding of the lower-level dual cannot cause infeasibility of \eqref{medial_level_dual_disc}.
Consequently, \cref{alg:ver} has the same convergence properties as the original algorithm.
According to \cite{jungenAdaptiveDiscretizationbasedAlgorithms2022a}, provided that the feasibility tolerance $\epsilon_{f}$ is greater than $0$, the algorithm is guaranteed to converge in a finite number of iterations under the assumptions that the host sets are compact, the defining functions are continuous, and the subproblems are solved to approximate global optimality.

In each iteration of \cref{alg:ver}, the dual formulation corresponding to an additional integer discretization point is added to problem \eqref{medial_level_dual_disc}.
Furthermore, the terms $\boldsymbol{\lambda}_{k}^T \mathbf{b}(\mathbf{y})$ and $\boldsymbol{\lambda}_{k}^T \mathbf{A}_{c}(\mathbf{y})$ lead to bilinear constraints with respect to the medial-level variables $\mathbf{y}$ and the dual variables $\boldsymbol{\lambda}_{k}$.
As a result, the optimization problem grows rapidly and becomes computationally challenging.
Consequently, the size of the problem that can be solved is limited, and we restrict the number of time steps per day in order to maintain tractability.

\subsection{Robustness of the identified design}
We verify the robustness of the candidate design for the illustrative case study using the method presented in the previous subsection.
Despite the removal of the electrical grid connection and modified parameters for the piecewise-linear input–output relations, the resulting system design is similar to the low total annualized cost design reported in the original publication by \cite{sassModelCompendiumData2020}.
The installed conversion components are a boiler, a small CHP plant, a heat pump, and a compression chiller.
The only installed storage component is a thermal energy storage (TES) unit for cooling.

Compared to the low-TAC design in \cite{sassModelCompendiumData2020}, the design identified in our study includes a heat pump but no heat storage.
This difference is likely due to the removal of the electrical grid connection.
Excess electricity produced by the CHP can no longer be exported to the grid and is instead converted into heat using the heat pump.
The installed capacities are given in Table \ref{tab:capacities}.
\begin{table}
    \begin{center}
        \begin{tabular}{l|c|c}
            \hline
            Component & Initial SOC fixed to \SI{50}{\%} of nominal capacity & Variable initial SOC\\
            \hline
            Boiler & \SI{441}{\kilo\watt} & \SI{433}{\kilo\watt}\\
            CHP & \SI{328}{\kilo\watt} & \SI{331}{\kilo\watt}\\
            Heat pump & \SI{153}{\kilo\watt} & \SI{150}{\kilo\watt}\\
            Compression Chiller & \SI{400}{\kilo\watt} & \SI{400}{\kilo\watt}\\
            TES cooling & \SI{1735}{\kilo\watt\hour} & \SI{1636}{\kilo\watt\hour}\\
            \hline
        \end{tabular}
        \caption{Nominal capacities of the installed components for the case where the initial state of charge (SOC) of the storage systems is fixed to \SI{50}{\%} of the nominal capacity and for the case where the initial SOC is a design variable.}
        \label{tab:capacities}
    \end{center}
\end{table}

Version $13.0.0$ of the Gurobi optimizer \citep{gurobi} is used to solve all subproblems on a $16$-core/$16$-thread Intel Xeon Gold 6248R CPU with \SI{3.0}{\giga\hertz}/\SI{4.0}{\giga\hertz} base/turbo frequency and \SI{64}{\giga\byte} of RAM running Linux Ubuntu version $24.04.2$ LTS.
For the design problem in the feasibility time-step heuristic, we set the `MIPGap' parameter to $0.005$, i.e., a relative optimality gap of \SI{0.5}{\%}, as solving the design problem to global optimality took prohibitively long, likely due to the large number of integer variables.
This optimality gap only affects the costs of the identified candidate design, the robustness of the design was verified with the default optimality gap of \SI{0.01}{\%}.
For the adaptive discretization algorithm (Section \ref{sec:method}), the feasibility tolerance, i.e., the value of the upper bound for which the design is determined robust, was chosen as $0.1$, and the optimality tolerance was chosen as $0.05$.
In total, seven iterations of the adaptive discretization algorithm (Section \ref{sec:method}) were required until the algorithm terminated due to the upper and lower bounds converging, i.e., seven realizations of the binary on/off variables were sufficient to cover all relevant uncertainty realizations.
This comparatively small number of discretization points is notable.
With four installed components and twelve time steps, the number of possible on/off states is $2^{4 \cdot 12}$.
The limited number of active discretization points is likely due to two factors.
First, the boiler and CHP are components with an unconstrained externally supplied input energy carrier and can therefore always be operated at maximum capacity.
Second, the installed cooling storage allows demand to be shifted between time steps, which reduces the number of distinct operational patterns required to cover the uncertainty realizations.

The solution of \eqref{medial_level_dual_disc} in the final iteration required \SI{721.3}{\second}, as reported by Gurobi \citep{gurobi}, and was solved to global optimality with an optimality gap of \SI{0.0}{\%}.
The largest optimality gap among all solutions of \eqref{medial_level_dual_disc} was \SI{0.0054}{\%}.
In total, solving \eqref{medial_level} required \SI{38959.0}{\second}.
Of this time, \SI{38868.3}{\second} were spent solving \eqref{medial_level_dual_disc}, whereas algorithm overhead and the solution of \eqref{lower_level} accounted for only \SI{90.7}{\second}.
Consequently, the solution of \eqref{medial_level_dual_disc} is the limiting factor regarding computational tractability.
With a growing number of discretization points, the required time to solve \eqref{medial_level_dual_disc} increases.
Overall, solution times for \eqref{medial_level_dual_disc} ranged from \SI{4.8}{\second} to \SI{18544.7}{\second}.
The first two iterations were solved substantially faster than the remaining ones, taking \SI{5.5}{\second} and \SI{4.8}{\second}, respectively, whereas all subsequent iterations required more than \SI{721.3}{\second}.
The number of required discretization points likely increases with problem size, as both an increased number of components and time steps lead to a larger number of combinations of possible on/off states.

Importantly, solving \eqref{medial_level} identifies a scenario within the convex hull of the historical data for which a violation of \SI{1.1}{\kilo\watt} occurs, i.e., the design is not robust.
The violation is small and occurs in the electricity balance; the demand at the corresponding time step is \SI{313.6}{\kilo\watt}.
To assess the sensitivity of the designs obtained using the feasibility time-step heuristic, we additionally determine designs based on $6$ and $8$ representative scenarios.
In both cases, the resulting designs are robust, with a maximum energy gap of \SI{0}{\kilo\watt}.

To model the energy storage, we assume a cyclical constraint, i.e., the last time step of each scenario is connected to the first.
Furthermore, we assume that in the first time step, the initial state of charge (SOC) is at \SI{50}{\%}.
This introduces an additional cooling demand even on days on which no cooling is required, as the energy dissipated by the TES needs to be compensated.
To investigate the influence of this assumption, we determine an additional design candidate using $4$ representative scenarios, treating the initial storage level as a long-term design variable in the upper-level problem. This formulation allows the optimizer to determine the optimal SOC at the beginning of each day. The resulting capacities are reported in the rightmost column of Table~\ref{tab:capacities}.

Overall, the design is very similar to the design without a variable initial state of charge, with slight deviations in capacity for most components.
The initial SOC of the TES for cooling is \SI{278}{\kilo\watt\hour}, which corresponds to \SI{17}{\%} and is less than in the fixed SOC case.
No infeasible scenarios were identified within the convex hull of the historical data by solving \eqref{medial_level} for this design, hence it is robust.

\section{Conclusion}\label{sec:conclusion}
Identifying energy system designs that can meet demand under any uncertainty realization in the predefined uncertainty set is crucial to ensure a reliable energy supply.
To this end, heuristic approaches such as the feasibility time-step heuristic \citep{teichgraeberExtremeEventsTime2020} have been applied, e.g., \cite{bahlRigorousSynthesisEnergy2017}.
This is because rigorous solution approaches, such as our robust energy system design approach \citep{wedemeyerRobustEnergySystem2025}, quickly become computationally intractable as system size and temporal resolution increase.
However, if the operational problem contains nonconvexities, the robustness of designs identified by this heuristic is not guaranteed.

This study provides a critical theoretical and empirical validation of the feasibility time-step heuristic in the context of MILP multi-energy system models with certain types of nonconvexities. Specifically, we identify three common nonconvexities in MILP multi-energy system models that affect robustness: piecewise-linear modeling of energy inflow-outflow behavior, minimum part-loads, and storage complementarity constraints.
Additionally, care must be taken to identify how the uncertain variables $\mathbf{y}$ enter the lower-level problem.
If the objective function of the lower-level problem is not jointly convex in $\mathbf{y}$ and $\mathbf{z}$, the feasibility time-step heuristic may fail even if the lower-level problem is convex in $\mathbf{z}$.

Our analysis shows that if curtailment of surplus energy flows is not allowed, any of these nonconvexities may lead to infeasible designs when using the feasibility time-step heuristic \citep{teichgraeberExtremeEventsTime2020}.
If curtailment is permitted, however, most of these issues disappear.
In particular, convex piecewise-linear input–output relations constitute only an apparent nonconvexity and can be reformulated as linear constraints, while storage complementarity can be ignored.
Under these conditions, minimum part-load constraints remain the primary source of potential infeasibility, although storage can partially mitigate their impact by shifting energy between time steps.

To rigorously verify the robustness of designs identified by the feasibility time-step heuristic, we propose a hybrid approach that combines adaptive discretization \citep{blankenshipInfinitelyConstrainedOptimization1976} with duality-based single-level reformulations \citep{fortuny-amatRepresentationEconomicInterpretation1981, burgardOptknockBilevelProgramming2003}.
While this approach provides a formal verification procedure, it leads to a quadratically constrained optimization problem whose size grows with each discretization iteration and is therefore computationally demanding.

We apply the hybrid approach to an illustrative case study based on \cite{sassModelCompendiumData2020} and \cite{vollAutomatedSuperstructurebasedSynthesis2013}.
A candidate design is first identified using the feasibility time-step heuristic \citep{teichgraeberExtremeEventsTime2020} and subsequently verified using the hybrid method.
A scenario within the convex hull of the historical data is identified for which a small infeasibility of \SI{1.1}{\kilo\watt} arises.
The mismatch occurs in the electricity balance, and the demand at the time step at which the violation occurs is \SI{313.6}{\kilo\watt}.
We additionally identify optimal designs with $6$ and $8$ representative days and find these designs to be robust.
Furthermore, we resolve the problem with $4$ representative scenarios while allowing the initial state of charge of storage components to be determined by the optimizer as a design variable, thereby identifying the optimal state of charge at the beginning of every day.
Despite the presence of minimum part-load constraints and the objective function being nonconvex in the nominal efficiency of the heat pump, which is uncertain due to its dependence on ambient temperature, the identified design is shown to be robust for all uncertainty realizations within the convex hull of the historical data.
Together with the theoretical analysis of the nonconvexities, these results demonstrate that the feasibility time-step heuristic can yield robust designs even when the operational problem is formulated as an MILP problem and nonconvexity in the uncertain variables is present.
However, the fact that the design for $4$ representative scenarios with a fixed initial state of charge was proven not to be robust highlights that care must be taken when applying the feasibility time-step heuristic.

While the present case study demonstrates that the feasibility time-step heuristic may identify non-robust designs of local multi-energy systems, future work may investigate larger systems to assess how frequently such failure modes arise in more complex applications.
However, the computational effort required for rigorous robustness verification increases with the number of possible combinations of on/off operating states.
Therefore, more efficient algorithms are needed before the proposed verification approach can be applied to substantially larger systems.

Future work should also investigate sources of nonconvexity that are beyond the scope of the present study.
In particular, nonconvexities arising from operational coupling constraints, i.e., limits on the number of start-ups and shut-downs, minimum up- and down-time, and ramping constraints, may affect the robustness of designs obtained with the feasibility time-step heuristic.

Moreover, we were unable to identify conversion components with inherently nonconvex input-output energy flow relations.
Future work may investigate whether there is a thermodynamic basis supporting this observation.

\section*{Acknowledgment}
This work was performed as part of the Helmholtz School for Data Science in Life, Earth and Energy (HDS-LEE) and received funding from the Helmholtz Association of German Research Centres. We further acknowledge financial support by the Helmholtz Association of German Research Centres through program-oriented funding.
Special thanks to Susanne Sass and Dinah Elena Hollermann for clarifying the mapping of the linearization parameters.

\section*{Data availability}
Data will be made available on request.

\section*{Declaration of competing interests}
The authors declare that they have no known competing financial interests or personal relationships that could have appeared to influence the work reported in this paper.

\section*{Authors' contributions}
Conceptualization: M.W., A.M., M.D.; Methodology: M.W.; Software: M.W.; Formal analysis and investigation: M.W.; Visualization: M.W.; Writing - original draft preparation: M.W.; Writing - review and editing: A.M., M.D.; Funding acquisition: A.M., M.D.; Supervision: A.M., M.D.

\section*{Declaration of generative AI and AI-assisted technologies in the manuscript preparation process}
During the preparation of this work, M.W. used ChatGPT to correct grammar and spelling and to improve the style of writing. After using this tool, all authors reviewed and edited the content as needed and take full responsibility for the content of the publication.
 
\section*{Nomenclature}
\subsection*{Abbreviations}
\begin{table}[H]
\begin{tabular}{ll}
    ESIP & Existence-constrained semi-infinite program\\
    GSIP & Generalized semi-infinite program\\
    MILP & Mixed-integer linear programming\\
    SOC & State of charge\\
    TAC & Total annualized cost\\
    TES & Thermal energy storage\\
\end{tabular}
\end{table}

\subsection*{Greek symbols}
\begin{table}[H]
\begin{tabular}{ll}
    $\beta$ & Slope parameter\\
    $\Delta_{t}$ & Time-step length\\
    $\epsilon_{f}$ & Feasibility tolerance\\
    $\epsilon_{o}$ & Optimality tolerance\\
    $\eta$ & Efficiency\\
    $\lambda$ & Linearization parameter\\
    $\boldsymbol{\lambda}$ & Dual variables\\
    $\tau$ & Time constant for energy loss\\
    $\phi$ & Auxiliary variable\\
\end{tabular}
\end{table}

\subsection*{Latin symbols}
\begin{table}[H]
\begin{tabular}{ll}
    $\mathbf{A}_{c}(\mathbf{y})$ & Constraint matrix for continuous variables\\
    $\mathbf{A}_{d}$ & Constraint matrix for integer variables\\
    $b$ & Binary variable\\
    $\mathbf{b}(\mathbf{y})$ & Right-hand side vector\\
    $c$ & Component\\
    $C$ & Set of components\\
    $\mathbf{c}_{c}$ & Cost vector\\
    $c_{inv}$ & Investment costs\\
    $c_{op}$ & Operational costs\\
    $d$ & Data point\\
    $\mathcal{D}$ & Historical data\\
    $e$ & Energy form\\
    $E$ & Energy\\
    $\mathcal{E}$ & Set of energy forms\\
    $\Edot$ & Energy flow\\
    $f_{lb}$ & Lower bound\\
    $f_{ub}$ & Upper bound\\
    $g$ & Constraint\\
    $g_{en}(\mathbf{x},\cdot, \cdot)$ & Constraints for energy system model\\
    $g_{x}(\mathbf{x})$ & Design constraints\\
    $g_{y}(\mathbf{y})$ & Uncertainty bounds\\
    $\mathcal{J}$ & Index set\\
    $k$ & Discretization index\\
    $\mathcal{K}$ & Set of discretization indices\\
    $n_{c}$ & Dimensionality of continuous operational variables\\
    \end{tabular}
    \end{table}
    \begin{table}[H]
    \begin{tabular}{ll}
    $n_{cs}$ & Number of constraints\\
    $n_{d}$ & Dimensionality of integer operational variables\\
    $n_{x}$ & Dimensionality of design variables\\
    $n_{y}$ & Dimensionality of uncertainty realizations\\
    $n_{z}$ & Dimensionality of operational variables\\
    $s$ & Scenario\\
    $\mathcal{S}$ & Set of scenarios\\
    $t$ & Time step\\
    $\mathcal{T}$ & Set of time steps\\
    $T_{amb}$ & Ambient temperature\\
    $\mathbf{x}$ & Design variables\\
    $\mathcal{X}$ & Feasible set of design variables\\
    $\mathbf{y}$ & Uncertainty realizations\\
    $\mathcal{Y}$ & Feasible set of uncertainty realizations\\
    $\mathbf{z}$ & Operational variables\\
    $\mathbf{z}_{c}$ & Continuous operational variables\\
    $\mathbf{z}_{d}$ & Integer operational variables\\
    $\mathcal{Z}$ & Feasible set of operational variables\\
    $\mathcal{Z}_{c}$ & Feasible set of operational variables for fixed integer variables\\
    \end{tabular}
\end{table}

\subsection*{Subscripts}
\begin{table}[H]
\begin{tabular}{ll}
    $agg$ & Aggregated\\
    $cur$ & Curtailment\\
    $dem$ & Demand\\
    $in$ & Inflow\\
    $hp$ & Heat pump\\
    $j$ & Index\\
    $k$ & Discretization index\\
    $max$ & Maximum\\
    $min$ & Minimum\\
    $nom$ & Nominal\\
\end{tabular}
\end{table}
\begin{table}[H]
\begin{tabular}{ll}
    $out$ & Outflow\\
    $s$ & Representative scenario\\
    $sto$ & Storage\\
    $t$ & Time step\\
\end{tabular}
\end{table}
 
\bibliographystyle{elsarticle-harv}
  \renewcommand{\refname}{References}  
  \bibliography{bibliography.bib}
\end{document}


\thispagestyle{firststyle}

  \begin{center}
    \begin{large}
      \textbf{\mytitle}
    \end{large} \\
    \myauthor
  \end{center}

  \vspace{0.5cm}

  \begin{footnotesize}
    \affil
  \end{footnotesize}

  \vspace{0.5cm}

\section{Parameters for the illustrative examples}
This section contains the parameters required to reproduce the examples from Sections $2.2$ and $2.3$ of the main manuscript.
\begin{table}[h]
    \begin{center}
        \begin{tabular}{c|c|c}
            \hline
            j & $\lambda_{in, j}$ & $\lambda_{out, j}$\\
            \hline
            0 & 0.1 & 0.075\\
            1 & 0.5 & 0.15\\
            2 & 1 & 0.95\\
            \hline
        \end{tabular}
        \caption{Parameters for the piecewise-linear nonconvex performance curve for the example from Section $2.2$ of the main manuscript.}
        \label{tab:lin_par_ex}
    \end{center}
\end{table}
\begin{table}[h]
    \begin{center}
        \begin{tabular}{c|c|c}
            \hline
            d & $\Edot_{1, dem, d}$ [\si{\kilo\watt}] & $\Edot_{2, dem, d}$ [\si{\kilo\watt}]\\
            \hline
            0 & 50 & 50\\
            1 & 80 & 3.75\\
            2 & 130 & 5\\
            \hline
        \end{tabular}
        \caption{Demand parameters for the example from Section $2.2$ of the main manuscript.}
        \label{tab:ex_par}
    \end{center}
\end{table}
\begin{table}[h]
    \begin{center}
        \begin{tabular}{c|c|c}
            \hline
            d & $\Edot_{1, dem, d}$ [\si{\kilo\watt}] & $\Edot_{2, dem, d}$ [\si{\kilo\watt}]\\
            \hline
            0 & 60 & 100\\
            1 & 150 & 8\\
            2 & 170 & 4\\
            \hline
        \end{tabular}
        \caption{Demand parameters for the example from Section $2.3$ of the main manuscript.}
        \label{tab:ex_par_2}
    \end{center}
\end{table}

\section{Linearization}\label{sec:lin}
\begin{figure}
    \centering
    \includegraphics[width=0.95\linewidth]{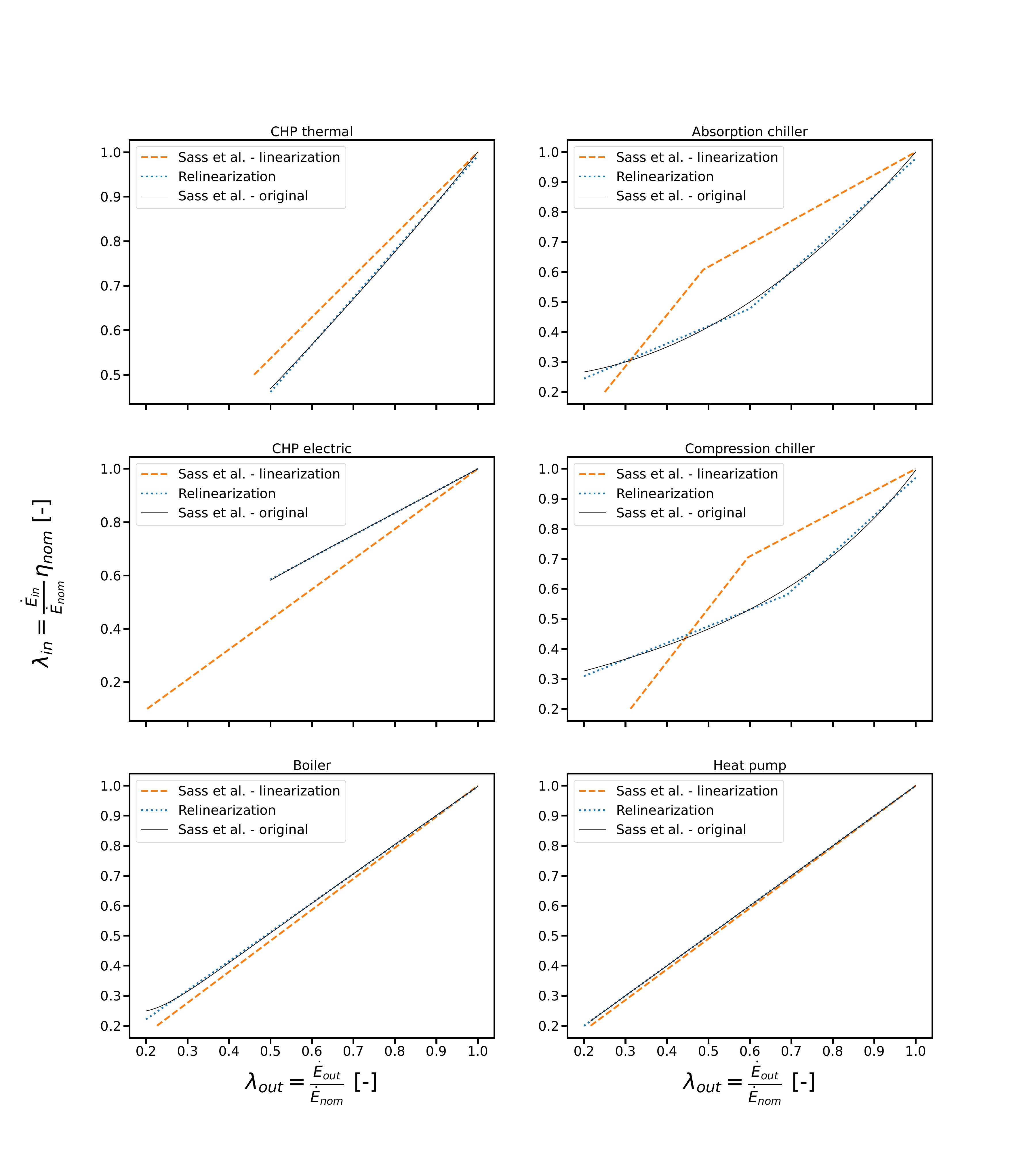}
    \caption{Original input-output relationship according to \cite{sassModelCompendiumData2020}, linearization by \cite{sassModelCompendiumData2020}, and our linearization (relinearization).}\label{fig:lin}
\end{figure}

As mentioned in the main manuscript, we relinearize the energy inflow-outflow for all conversion components, i.e., boiler, CHP, absorption chiller (AC), compression chiller (CC), and the heat pump (HP), in the case study by \cite{sassModelCompendiumData2020}.
To ensure that the piecewise energy inflow-outflow relationship closely approximates the original nonlinear relationship, we apply the following linearization procedure:
We sample the nonlinear function at $250$ equidistant points between the minimum part-load and the maximum part-load.
We then minimize the sum of the squared residuals between a piecewise-linear function and the sampled points.
The locations of the breakpoints of the piecewise-linear function are optimization variables, and we fix the output part-load of the lowest and highest breakpoints to be at the minimum and maximum part-load, respectively.
To ensure comparability, we use the same number of breakpoints as \cite{sassModelCompendiumData2020}.
Table \ref{tab:lin_par} shows the parameters of the linearization, i.e., the relative input $\lambda_{in}$, which is the energy inflow $\Edot_{in}$ divided by the nominal capacity $\Edot_{nom}$ multiplied by the nominal efficiency $\eta_{nom}$, and the relative output $\lambda_{out}$ which is calculated from the energy outflow $\Edot_{out}$ and the nominal capacity.
Figure \ref{fig:lin} shows the original input-output relationship, the linearization by \cite{sassModelCompendiumData2020}, and our linearization (relinearization).

\begin{table}[h]
    \begin{center}
        \begin{tabular}{c|c|c|c|c|c|c}
            \hline
            Component & $\lambda_{in, 0}$ & $\lambda_{out, 0}$ & $\lambda_{in, 1}$ & $\lambda_{out, 1}$ & $\lambda_{in, 2}$ & $\lambda_{out, 2}$ \\
            \hline
            Boiler & 0.221 & 0.2 & 0.996 & 1 & - & -\\
            CHP thermal & 0.585 & 0.5 & 1 & 1 & - & -\\
            CHP electric & 0.461 & 0.5 & 0.992 & 1 & - & -\\
            Absorption chiller (AC) & 0.245 & 0.2 & 0.478 & 0.6 & 0.978 & 1\\
            Compression chiller (CC) & 0.309 & 0.2 & 0.580 & 0.689 & 0.97 & 1\\
            Heat pump (HP) & 0.2 & 0.2 & 1 & 1 & - & -\\
            \hline
        \end{tabular}
        \caption{Linearization parameters for all conversion components.}
        \label{tab:lin_par}
    \end{center}
\end{table}
\section{Optimization models}
In general, all components are modeled according to the description by \cite{sassModelCompendiumData2020}.
The relevant differences between our study and \cite{sassModelCompendiumData2020} are:
(i) the removal of the electricity grid connection, (ii) our parameters for the linearized performance curves (cf. Section \ref{sec:lin}), and (iii) the utilized representative scenario data and extreme scenarios, which are determined according to the description in Section $3$ of the main manuscript.

Throughout the manuscript, scalar-valued quantities are denoted in regular font, e.g., $x$, vector-valued quantities are denoted in bold font, e.g., $\mathbf{x}$, and set-valued quantities are denoted in calligraphic font, e.g., $\mathcal{X}$.

\begingroup
\begin{align*}
    \underset{\mathbf{x}, \mathbf{z}}{\min} & \sum_{c \in \mathcal{C}}\left(Capex_{c} \left( \frac{1}{\gamma_{pvf}} + \gamma_{maintenance, c} \right)\right) +  \sum_{c \in \mathcal{C}_{fuel}} \sum_{s \in \mathcal{S}} \omega_{s} \sum_{t \in \mathcal{T}}\left(\gamma_{fuel} \Edot_{fuel, c, s, t} \Delta_{t}\right)&\\
    \text{s.t.} \quad & \Edot_{e, dem, s, t} - \sum_{c \in \mathcal{C}_{e, out}} \Edot_{e, c, s, t} + \sum_{c \in \mathcal{C}_{e, in}} \Edot_{e, c, s, t} + \sum_{c \in \mathcal{C}_{e, sto}} \left(\Edot_{e, c, s, t, in} - \Edot_{e, c, s, t, out} \right) \le 0 \ & \forall e \in \mathcal{E}, \forall s \in \mathcal{S}, \forall t \in \mathcal{T}\\
    & \Edot_{nom, c} b_{ex, c} \le \Edot_{max, c} & \forall c \in \mathcal{C} \\
    & \Edot_{nom, c} b_{ex, c} \ge \Edot_{min, c} & \forall c \in \mathcal{C} \\
    & Capex_{c} = \sum_{j \in \mathcal{J}_{capex, c}} b_{capex, c, j} Capex_{c, lb, j} + \beta_{capex, c, j} \left(\Edot_{nom, c, j} - \Edot_{lb, c, j} b_{capex, c, j}\right) & \forall c \in \mathcal{C}\\
    & \Edot_{nom, c} = \sum_{j \in \mathcal{J}_{capex, c}} \Edot_{nom, c, j} & \forall c \in \mathcal{C}\\
    & \Edot_{nom, c, j} b_{capex, c, j} \le \Edot_{lb, c, j + 1} & \forall c \in \mathcal{C}, \forall j \in \mathcal{J}_{capex, c} \\
    & \Edot_{nom, c, j} b_{capex, c, j} \ge \Edot_{lb, c, j} & \forall c \in \mathcal{C}, \forall j \in \mathcal{J}_{capex, c}\\
    & \sum_{j \in \mathcal{J}_{capex, c}} b_{capex, c, j} \le 1 & \forall c \in \mathcal{C}\\
    & \Edot_{e_{in}(c), c, s, t} = \sum_{j \in \mathcal{J}_{eff, c}} b_{eff, e, c, s, t, j} \lambda_{in, e, c, j} \frac{\Edot_{nom, c}}{\eta_{nom, c, s, t}}\\
    & + \frac{\beta_{eff, e, c, j}}{\eta_{nom, c, s, t}}\left(\Edot_{e, c, s, t, j} - \lambda_{out, e, c, j}  b_{eff, e, c, s, t, j} \Edot_{nom, c}\right) & \mathllap{\forall c \in \mathcal{C}_{conv}, \forall e \in \mathcal{E}_{out}(c), \forall s \in \mathcal{S}, \forall t \in \mathcal{T}}\\
    & \Edot_{e, c, s, t} = \sum_{j \in \mathcal{J}_{eff, c}} \Edot_{e, c, s, t, j} & \mathllap{\forall c \in \mathcal{C}_{conv}, \forall e \in \mathcal{E}_{out}(c), \forall s \in \mathcal{S}, \forall t \in \mathcal{T}}\\
    & \lambda_{out, e, c, j} \Edot_{nom, c} b_{eff, e, c, s, t, j} \le \Edot_{e, c, s, t, j} & \mathllap{\forall c \in \mathcal{C}_{conv}, \forall j \in \mathcal{J}_{eff, c}, \forall e \in \mathcal{E}_{out}(c), \forall s \in \mathcal{S}, \forall t \in \mathcal{T}}\\
    & \lambda_{out, e, c, j + 1} \Edot_{nom, c} b_{eff, e, c, s, t, j} \ge \Edot_{e, c, s, t, j} & \mathllap{\forall c \in \mathcal{C}_{conv}, \forall j \in \mathcal{J}_{eff, c}, \forall e \in \mathcal{E}_{out}(c), \forall s \in \mathcal{S}, \forall t \in \mathcal{T}}\\
    & \sum_{j \in \mathcal{J}_{eff, c}} b_{eff, e, c, s, t, j} \le 1 & \mathllap{\forall c \in \mathcal{C}_{conv}, \forall e \in \mathcal{E}_{out}(c), \forall s \in \mathcal{S}, \forall t \in \mathcal{T}}\\
    & E_{e(c), c, s, t} \le \gamma_{sto, c}\Edot_{nom, c} & \forall c \in \mathcal{C}_{sto}, \forall s \in \mathcal{S}, \forall t \in \mathcal{T}\\
    & \Edot_{e(c), c, s, t, in} \le \Edot_{nom, c} & \forall c \in \mathcal{C}_{sto}, \forall s \in \mathcal{S}, \forall t \in \mathcal{T}\\
    & \Edot_{e(c), c, s, t, out} \le \Edot_{nom, c} & \forall c \in \mathcal{C}_{sto}, \forall s \in \mathcal{S}, \forall t \in \mathcal{T}\\
    & E_{e(c), c, s, t} \left( 1 + \frac{\Delta_{t}}{\tau_{c}}\right) = E_{e(c), c, s, t - 1} + \Delta_{t} \left(\eta_{c, in} \Edot_{e(c), c, s, t, in} - \frac{\Edot_{e(c), c, s, t, out}}{\eta_{c, out}} \right) & \mathllap{\forall c \in \mathcal{C}_{sto}, \forall s \in \mathcal{S}, \forall t \in \mathcal{T} \setminus \{1\}}\\
    & E_{e(c), c, s, 1} \left( 1 + \frac{\Delta_{t}}{\tau_{c}}\right) = E_{e(c), c, s, |\mathcal{T}|} + \Delta_{t} \left(\eta_{c, in} \Edot_{e(c), c, s, 1, in} - \frac{\Edot_{e(c), c, s, 1, out}}{\eta_{c, out}} \right) & \mathllap{\forall c \in \mathcal{C}_{sto}, \forall s \in \mathcal{S}}\\
     & E_{e(c), c, s, 1} = E_{init, c} & \forall c \in \mathcal{C}_{sto}, \forall s \in \mathcal{S}\\
      & \Edot_{e(c), c, s, t} \le f_{solar, s, t} \Edot_{nom, c} & \mathllap{\forall c \in \left\{PV_{i} \ | \ \forall i \in 1, \dots, n_{components}\right\}, \forall s \in \mathcal{S}, \forall t \in \mathcal{T}}\\
\end{align*}
\endgroup

with
\begin{align*}
    &\mathbf{x} = [\underbrace{Capex_{c}, \Edot_{nom, c}, b_{ex, c}}_{\forall c \in \mathcal{C}}, \underbrace{E_{init, c}}_{\forall c \in \mathcal{C}_{sto}},    \underbrace{\Edot_{nom, c, j}, b_{capex, c, j}}_{\forall c \in \mathcal{C}, \forall j \in \mathcal{J}_{capex, c}}] \\
    & Capex_{c} \in \mathbb{R}_{\ge 0}, \Edot_{nom, c} \in \mathbb{R}_{\ge 0}, b_{ex, c} \in \{0, 1\} \\
    &\mathbf{z} = [\underbrace{\Edot_{e, c, s, t}}_{\forall c \in \mathcal{C}_{conv}, \forall e \in \mathcal{E}(c), \forall s \in \mathcal{S}, \forall t \in \mathcal{T}}, \underbrace{\Edot_{e, c, s, t, in}, \Edot_{e, c, s, t, out}}_{\forall c \in \mathcal{C}_{sto}, \forall e \in \mathcal{E}(c), \forall s \in \mathcal{S}, \forall t \in \mathcal{T}}, \underbrace{E_{e(c), c, s, t}}_{\forall c \in \mathcal{C}_{sto}, \forall s \in \mathcal{S}, \forall t \in \mathcal{T}}, \underbrace{\Edot_{e, c, s, t, j}, b_{eff, e, c, s, t, j}}_{\forall c \in \mathcal{C}_{conv}, \forall e \in \mathcal{E}_{out}(c), \forall s \in \mathcal{S}, \forall t \in \mathcal{T}, \forall j \in \mathcal{J}_{eff, c}}] \\
    & \Edot_{e, c, s, t} \in \mathbb{R}_{\ge 0},
    \Edot_{e, c, s, t, j} \in \mathbb{R}_{\ge 0},\Edot_{e, c, s, t, in} \in \mathbb{R}_{\ge 0}, \Edot_{e, c, s, t, out} \in \mathbb{R}_{\ge 0},
    E_{e(c), c, s, t} \in \mathbb{R}_{\ge 0}, b_{eff, e, c, s, t, j} \in \{0, 1\}\\
    & \mathcal{C} = \{CHP_{small, i}, CHP_{medium, i}, CHP_{large, i}, Boiler_{i}, CC_{i}, AC_{i}, HP_{i}, PV_{i}, Battery_{i},\\
    & \quad TES_{heat, i}, TES_{cooling, i} \ | \ i = 1, \dots, n_{components} \}\\
    & \mathcal{C}_{fuel} = \left\{CHP_{small, i}, CHP_{medium, i}, CHP_{large, i}, Boiler_{i} \ | \ i = 1, \dots, n_{components} \right\}\\
    & \mathcal{C}_{conv} = \left\{CHP_{small, i}, CHP_{medium, i}, CHP_{large, i}, Boiler_{i}, CC_{i}, AC_{i}, HP_{i}\ | \ i = 1, \dots, n_{components} \right\}\\
    & \mathcal{C}_{sto} = \left\{Battery_{i}, TES_{heat, i}, TES_{cooling, i} \ | \ i = 1, \dots, n_{components} \right\}\\
    & \mathcal{C}_{heat, out} = \left\{CHP_{small, i}, CHP_{medium, i}, CHP_{large, i}, Boiler_{i}, HP_{i} \ | \ i = 1, \dots, n_{components} \right\}\\
    & \mathcal{C}_{heat, in} = \left\{AC_{i} \ | \ i = 1, \dots, n_{components} \right\}\\
    & \mathcal{C}_{electricity, out} = \left\{CHP_{small, i}, CHP_{medium, i}, CHP_{large, i} \ | \ i = 1, \dots, n_{components} \right\}\\
    & \mathcal{C}_{electricity, in} = \left\{CC_{i}, HP_{i} \ | \ i = 1, \dots, n_{components} \right\}\\
    & \mathcal{C}_{cooling, out} = \left\{CC_{i}, AC_{i} \ | \ i = 1, \dots, n_{components} \right\}\\
    & \mathcal{E} = \left\{heat, electricity, cooling \right\}\\
\end{align*}
Here, $\mathbf{x}$ are the design variables and $\mathbf{z}$ are the operational variables.
The considered components are boilers, three different sizes of CHPs (small, medium, large), absorption chillers (ACs), compression chillers (CCs), heat pumps (HPs), photovoltaic modules (PVs), batteries, and thermal energy storage (TES) for heat and cooling.
$Capex_{c}$ are the capital expenditures to install component $c$, $\Edot_{fuel, c, s, t}$ are the energy imports of fuel, i.e., natural gas, for the fuel burning components, which are boilers and CHPs.
In our example, up to $n_{components} = 2$ copies can be installed of each component.

$\mathcal{E}$ is the set of energy forms, $\mathcal{S}$ is the set of scenarios and $\mathcal{T}$ is the set of time steps.
$\mathcal{J}_{capex, c}$ is the set of linearization indices for the linearization of the capital expenditures for each component, and $\mathcal{J}_{eff, c}$ is the set of linearization indices for the linearization of the energy inflow-outflow relationship for each component. 
$\mathcal{C}$ is the set of components, $\mathcal{C}_{fuel}$ is the set of components supplied externally by gas, $\mathcal{C}_{conv}$ is the set of conversion components, $\mathcal{C}_{sto}$ is the set of storage components, $\mathcal{C}_{heat, out}$ and $\mathcal{C}_{heat, in}$ are the sets of components whose energy outflows and inflows are heat, respectively.
$\mathcal{C}_{electricity, out}$ and $\mathcal{C}_{electricity, in}$ are the sets of components whose energy outflows and inflows are electricity, respectively, and $\mathcal{C}_{cooling, out}$ is the set of components whose energy outflows are cooling and there are no conversion components taking cooling energy as input.
$e_{in}(c)$ is a function that maps each component to its inflow energy, and $e(c)$ maps each storage component to its inflow and outflow energy.
A similar mapping cannot be constructed for the energy outflows, as the CHPs have two output energy forms.

$\Edot_{e, c, s, t}$ are the energy flows of energy form $e$ for conversion component $c$ in scenario $s$ and time step $t$, while $\Edot_{e, c, s, t, in}$ and $\Edot_{e, c, s, t, out}$ are the energy inflows and outflows of storage components, and $E_{init, c}$ is the initial state of charge.
In the default case study, the initial state of charge is fixed to $50\%$ of the nominal capacity; in the variable SOC case study, it is a design variable.
$\Edot_{nom, c}$ is the nominal capacity of component $c$ and $b_{ex, c}$ is a binary variable that indicates whether component $c$ has been installed.

The nonlinear dependency of $Capex_{c}$ on  $\Edot_{nom, c}$ is approximated by a piecewise-linear relationship.
To this end, binary variables $b_{capex, c, j}$ indicate which segment $j \in \mathcal{J}_{capex, c}$ of the piecewise linearization is active.
$\Edot_{nom, c, j}$ is the nominal capacity variable for each respective segment.
Similarly, the nonlinear dependency of $\Edot_{e_{in}(c), c, s, t}$ on $\Edot_{e, c, s, t}$ is approximated by a piecewise-linear relationship.
Importantly, for the CHPs, there are two output energy forms, electricity and heat.
Hence, $\mathcal{E}_{out}(c)$ is the index set that contains all output energy forms for each component.
Again, $b_{eff, e, c, s, t, j}$ are binary variables that indicate which segment $j \in \mathcal{J}_{eff, c}$ of the linearization is active, and $\Edot_{e, c, s, t, j}$ is the output energy flow for that segment.
$E_{e(c), c, s, t}$ are the storage levels of each storage.
The equation to determine the energy storage level is obtained by using the implicit Euler scheme on the differential equation of the storage.
For a detailed description of the linearization formulation and the modeling equations, we refer to \cite{sassModelCompendiumData2020}.
Note that where \cite{sassModelCompendiumData2020} use $Q$ to refer to energy flow, we use $\Edot$.

The energy flow demands $\Edot_{e, dem, s, t}$, the capacity factor of the PV modules $f_{solar, s, t}$, and the ambient temperature $T_{amb, s, t}$ are uncertain variables that are determined for each representative scenario $s$ from historical data.
The ambient temperature is then used to calculate the efficiency of the heat pumps: $\eta_{nom, HP_{i}, s, t}(T_{amb, s, t}) = \frac{0.36T_{hp}}{T_{hp} - T_{amb, s, t}}$.
For all other components, the nominal efficiency $\eta_{nom, c, s, t}$ is constant and does not depend on uncertain variables.

The objective function minimizes the total annualized cost.
To this end, the present value factor is defined as
\begin{equation*}
    \gamma_{pvf} = \frac{\left(1 + 0.08\right)^{4} - 1}{0.08\left(1 + 0.08\right)^{4}}
\end{equation*}
for an interest rate of $0.08$ and a time horizon of $4$ years.
The maintenance factors $\gamma_{maintenance, c}$ can be obtained from Tables $1$ and $2$ in \cite{sassModelCompendiumData2020}.
The scenario weights are defined as
\begin{equation*}
    \omega_{s} = 365 \frac{n_{data, s}}{n_{data}},
\end{equation*}
where $n_{data, s}$ are the number of days in the historical data in the cluster of scenario $s$ and $n_{data}$ are the total number of days in the historical data.
$\gamma_{fuel} = \SI{0.08}{\frac{EUR}{\kilo\watt\hour}}$ are the fuel costs and $\Delta_{t}$ is the time-step length, which is \SI{2}{\hour} in our example.

The parameters for the economic variables $Capex_{c, lb, j}$ and the capacities $\Edot_{lb, c, j}$ can be obtained from \cite{sassModelCompendiumData2020} in their publication from Table C.2.
Note, that we use $\Edot_{lb, c, j}$ where \cite{sassModelCompendiumData2020} use $Q^{lb}_{i, j}$, $\Edot_{max, c}$ and $\Edot_{min, c}$ correspond to the smallest and largest values of $\Edot_{lb, c, j}$.
$\beta_{capex, c, j}$ is calculated as
\begin{equation*}
    \beta_{capex, c, j} = \frac{Capex_{c, lb, j + 1} - Capex_{c, lb, j}}{\Edot_{lb, c, j + 1} - \Edot_{lb, c, j}}.
\end{equation*}

The values for the linearization parameters $\lambda_{in, e, c, j}$ and $\lambda_{out, e, c, j}$ of the piecewise-linear energy inflow-outflow curves can be found in Section \ref{sec:lin}.
$\beta_{eff, e, c, j}$ is calculated as
\begin{equation*}
    \beta_{eff, e, c, j} = \frac{\lambda_{in, e, c, j + 1} - \lambda_{in, e, c, j}}{\lambda_{out, e, c, j + 1} - \lambda_{out, e, c, j}}.
\end{equation*}

$\gamma_{sto, c}$ is the coefficient that determines the relation between the nominal capacity of each storage component and the maximum storage capacity.
For batteries it is $\SI{4}{\frac{\kilo\watt\hour}{\kilo\watt}}$ and for thermal energy storages $\SI{1}{\frac{\kilo\watt\hour}{\kilo\watt}}$.
While $\tau_{c}$ are the time constants for self-discharge, and $\eta_{c, in}$ and $\eta_{c, out}$ are the charging and discharging efficiencies.

The operational problem is formulated as:
\begingroup
\begin{align*}
    \underset{\mathbf{z}, \phi}{\min} & \ \phi&\\
    \text{s.t.} \quad & \Edot_{e, dem, t} - \sum_{c \in \mathcal{C}_{e, out}} \Edot_{e, c, t} + \sum_{c \in \mathcal{C}_{e, in}} \Edot_{e, c, t} + \sum_{c \in \mathcal{C}_{e, sto}} \left(\Edot_{e, c, t, in} - \Edot_{e, c, t, out} \right) \le \phi \ & \forall e \in \mathcal{E}, \forall t \in \mathcal{T}\\
    & \Edot_{e_{in}(c), c, t} = \sum_{j \in \mathcal{J}_{eff, c}} b_{eff, e, c, t, j} \lambda_{in, e, c, j} \frac{\Edot_{nom, c}}{\eta_{nom, c, t}}\\
    & + \frac{\beta_{eff, e, c, j}}{\eta_{nom, c, t}}\left(\Edot_{e, c, t, j} - \lambda_{out, e, c, j}  b_{eff, e, c, t, j} \Edot_{nom, c}\right) & \mathllap{\forall c \in \mathcal{C}_{conv}, \forall e \in \mathcal{E}_{out}(c), \forall t \in \mathcal{T}}\\
    & \Edot_{e, c, t} = \sum_{j \in \mathcal{J}_{eff, c}} \Edot_{e, c, t, j} & \mathllap{\forall c \in \mathcal{C}_{conv}, \forall e \in \mathcal{E}_{out}(c), \forall t \in \mathcal{T}}\\
    & \lambda_{out, e, c, j} \Edot_{nom, c} b_{eff, e, c, t, j} \le \Edot_{e, c, t, j} & \mathllap{\forall c \in \mathcal{C}_{conv}, \forall j \in \mathcal{J}_{eff, c}, \forall e \in \mathcal{E}_{out}(c), \forall t \in \mathcal{T}}\\
    & \lambda_{out, e, c, j + 1} \Edot_{nom, c} b_{eff, e, c, t, j} \ge \Edot_{e, c, t, j} & \mathllap{\forall c \in \mathcal{C}_{conv}, \forall j \in \mathcal{J}_{eff, c}, \forall e \in \mathcal{E}_{out}(c), \forall t \in \mathcal{T}}\\
    & \sum_{j \in \mathcal{J}_{eff, c}} b_{eff, e, c, t, j} \le 1 & \mathllap{\forall c \in \mathcal{C}_{conv}, \forall e \in \mathcal{E}_{out}(c), \forall t \in \mathcal{T}}\\
    & E_{e(c), c, t} \le \gamma_{sto, c}\Edot_{nom, c} & \forall c \in \mathcal{C}_{sto}, \forall t \in \mathcal{T}\\
    & \Edot_{e(c), c, t, in} \le \Edot_{nom, c} & \forall c \in \mathcal{C}_{sto}, \forall t \in \mathcal{T}\\
    & \Edot_{e(c), c, t, out} \le \Edot_{nom, c} & \forall c \in \mathcal{C}_{sto}, \forall t \in \mathcal{T}\\
    & E_{e(c), c, t} \left( 1 + \frac{\Delta_{t}}{\tau_{c}}\right) = E_{e(c), c, t - 1} + \Delta_{t} \left(\eta_{c, in} \Edot_{e(c), c, t, in} - \frac{1}{\eta_{c, out}} \Edot_{e(c), c, t, out}\right) & \forall c \in \mathcal{C}_{sto}, \forall t \in \mathcal{T} \setminus \{1\}\\
    & E_{e(c), c, 1} \left( 1 + \frac{\Delta_{t}}{\tau_{c}}\right) = E_{e(c), c, |\mathcal{T}|} + \Delta_{t} \left(\eta_{c, in} \Edot_{e(c), c, 1, in} - \frac{1}{\eta_{c, out}} \Edot_{e(c), c, 1, out}\right) & \forall c \in \mathcal{C}_{sto}\\
    & E_{e(c), c, 1} = E_{init, c} & \forall c \in \mathcal{C}_{sto}\\
    & \Edot_{e(c), c, t} \le f_{solar, t} \Edot_{nom, c} & \mathllap{\forall c \in \left\{PV_{i} \ | \ \forall i \in 1, \dots, n_{components}\right\}, \forall s \in \mathcal{S}, \forall t \in \mathcal{T}}\\
\end{align*}
\endgroup

with
\begin{align*}
    &\mathbf{z} = [\underbrace{\Edot_{e, c, t}}_{\forall c \in \mathcal{C}_{conv}, \forall e \in \mathcal{E}(c), \forall t \in \mathcal{T}}, \underbrace{\Edot_{e, c, t, in}, \Edot_{e, c, t, out}}_{\forall c \in \mathcal{C}_{sto}, \forall e \in \mathcal{E}(c), \forall t \in \mathcal{T}}, \underbrace{E_{e(c), c, t}}_{\forall c \in \mathcal{C}_{sto}, \forall t \in \mathcal{T}}, \underbrace{\Edot_{e, c, t, j}, b_{eff, e, c, t, j}}_{\forall c \in \mathcal{C}_{conv}, \forall e \in \mathcal{E}_{out}(c), \forall t \in \mathcal{T}, \forall j \in \mathcal{J}_{eff, c}}] \\
    & \phi \in \mathbb{R}, \Edot_{e, c, t} \in \mathbb{R}_{\ge 0},
    \Edot_{e, c, t, j} \in \mathbb{R}_{\ge 0},\Edot_{e, c, t, in} \in \mathbb{R}_{\ge 0}, \Edot_{e, c, t, out} \in \mathbb{R}_{\ge 0},
    E_{e(c), c, t} \in \mathbb{R}_{\ge 0}, b_{eff, e, c, t, j} \in \{0, 1\}\\
\end{align*}

The operational problem is similar to the design problem.
The equations to calculate the capital expenditures are omitted, and the objective function is replaced by an auxiliary variable $\phi$, which measures the maximum violation of the energy balance equations.
Variables and parameters are similar to the design problem, except that no scenario index is included, as there is only one scenario per operational problem.
 
\bibliographystyle{apalike}
  \renewcommand{\refname}{References}  
  \bibliography{bibliography.bib}